\documentclass[12pt]{amsart}
\usepackage{amsmath,amsfonts,amssymb,amsthm, multicol}
\usepackage{anysize}
\marginsize{2cm}{2cm}{2cm}{2cm}
\usepackage{graphicx}
\usepackage{color}
\usepackage[pdftex]{hyperref}
\usepackage{enumerate}
\usepackage{caption}

\usepackage{mathtools}

\DeclarePairedDelimiter\floor{\lfloor}{\rfloor}

\def\g{\mathfrak{g}}

\def\n{\mathfrak{n}}

\def\R{\mathbb{R}}
\def\Q{\mathbb{Q}}
\def\Z{\mathbb{Z}}
\def\N{\mathbb{N}}

\def\e{\operatorname{e}}

\def\ad{\operatorname{ad}}
\def\tr{\operatorname{tr}}

\def\alt{\raise1pt\hbox{$\bigwedge$}}
\def\pint{\langle \cdotp,\cdotp \rangle }

\def\vol{\operatorname{vol}}

\def\mid{\, \vert \,} 
\def\Span{\operatorname{span}}
\def\diag{\operatorname{diag}}

\theoremstyle{plain}
\newtheorem{theorem}{\bf Theorem}[section]
\newtheorem{corollary}[theorem]{\bf Corollary}
\newtheorem{proposition}[theorem]{\bf Proposition}
\newtheorem{lemma}[theorem]{\bf Lemma}

\theoremstyle{definition}
\newtheorem{definition}[theorem]{\bf Definition}
\newtheorem{example}[theorem]{\bf Example}
\newtheorem{hypothesis}{\bf Hypothesis}
\newtheorem{hipotesis}{\bf Hypothesis}

\theoremstyle{remark}
\newtheorem{remark}[theorem]{\bf Remark}

\newcommand{\ri}{{\rm (i)}}
\newcommand{\rii}{{\rm (ii)}}
\newcommand{\riii}{{\rm (iii)}}
\newcommand{\riv}{{\rm (iv)}}

\title[Symplectic solvmanifolds satisfying the hard-Lefschetz condition]{Construction of symplectic solvmanifolds satisfying the hard-Lefschetz condition}

\author{Adrián Andrada}
\email{adrian.andrada@unc.edu.ar}
\author{Agustín Garrone}
\email{agustin.garrone@mi.unc.edu.ar}

\date{}
\address{FAMAF, Universidad Nacional de C\'ordoba and CIEM-CONICET, Av. Medina Allende s/n, Ciudad Universitaria, X5000HUA C\'ordoba, Argentina}

\subjclass[2020]{53D05, 22E25, 22E40}
\keywords{Almost Kähler solvmanifold, symplectic form, hard-Lefschetz condition, almost abelian Lie group, lattice}

\begin{document}

\begin{abstract} 
A compact symplectic manifold $(M, \omega)$ is said to satisfy the hard-Lefschetz condition if it is possible to develop an analogue of Hodge theory for $(M, \omega)$. This loosely means that there is a notion of harmonicity of differential forms in $M$, depending on $\omega$ alone, such that every de Rham cohomology class in has a $\omega$-harmonic representative. In this article, we study two non-equivalent families of diagonal almost-abelian Lie algebras that admit a distinguished almost-Kähler structure and compute their cohomology explicitly. We show that they satisfy the hard-Lefschetz condition with respect to any left-invariant symplectic structure by exploiting an unforeseen connection with Kneser graphs. We also show that for some choice of parameters their associated simply connected, completely solvable Lie groups admit lattices, thereby constructing examples of almost-Kähler solvmanifolds satisfying the hard-Lefschetz condition, in such a way that their de Rham cohomology is fully known.
\end{abstract}

\maketitle

\section{Introduction} 
\indent On any compact symplectic manifold $(M,\omega)$ of dimension $2n$, whose algebra of differential forms is denoted by $\Omega(M)$, it makes sense to consider the \textit{Lefschetz operators} $L_m^{(n)}:\Omega^m(M) \to \Omega^{2n-m}(M)$ given by $L_m^{(n)}(\alpha) := \omega^m \wedge \alpha$, where $0 \leq m \leq n$. They induce operators on the de Rham cohomology $L_m^{(n)}:H_{dR}^m(M) \to H_{dR}^{2n-m}(M)$, denoted by the same way by a slight abuse of language and defined simply by $L_m^{(n)}([\alpha]) = [ \omega^m \wedge \alpha]$. It is said that $(M, \omega)$ \textit{satisfies the hard-Lefschetz condition} if $L_m^{(n)}$ is a linear isomorphism for all $0 \leq m \leq n$. We elaborate on the importance of finding examples of compact symplectic manifolds that satisfy the hard-Lefschetz condition in Section \ref{seccion: the hard lefschetz condition}. We are particularly interested in the situation in which $M$ is a \textit{solvmanifold}, which is just a compact quotient $\Gamma \backslash G$ of some simply connected  solvable Lie group $G$ by a discrete co-compact subgroup $\Gamma$ of $G$, called a \textit{lattice}. A short review of them is given in Section \ref{seccion: solvmanifolds}. We mention here that a valuable tool in the study of the cohomology of solvmanifolds is \textit{Hattori's theorem}, which states that if $G$ is a completely solvable Lie group %the natural inclusion $\alt^* \g^* \hookrightarrow \Omega^*(\Gamma\backslash G)$ induces 
then there is an isomorphism in cohomology
$H^*(\g) \cong H^*_{dR}(\Gamma\backslash G)$, where $\g$ is the Lie algebra of $G$. We focus on some \textit{almost abelian} Lie groups $G$, defined as the simply connected Lie groups corresponding to a Lie algebra $\g$ that has a codimension one abelian ideal. Upon a choice of basis of $\g$, these can be described via the semidirect product $\R \ltimes_A \R^{2n-1}$ where $A \in \mathfrak{gl}(2n-1, \mathbb{R})$ is a matrix carrying the information of all nonzero brackets; we call $\g_A := \R \ltimes_A \R^{2n-1}$ this Lie algebra, and also $G_A := \mathbb{R} \ltimes_{\phi} \R^{2n-1}$ the corresponding simply connected Lie group to $\g_A$, where $\phi(t) := \exp(tA)$. More on them is said in the first paragraphs in Section \ref{seccion: almost abelian lie groups}. We are interested in diagonal matrices
\begin{equation}\label{eq: matrixA}
A=\diag(0, b_2, \ldots, b_n, - b_2, \ldots, -b_n)
\end{equation}
that satisfy either one of the following assumptions:
\begin{hipotesis} \label{hip: 1}
    No finite nontrivial addition or subtraction of different elements in the set $\{b_l \mid 2 \leq l \leq n\}$ is zero. That is, if $\sum_{j=2}^n \varepsilon_jb_j=0$ for some $\varepsilon_j \in \{-1,0,1\}$ then $\varepsilon_j=0$ for all $j$.
\end{hipotesis}
  %  \indent Note that Hypothesis \ref{hip: 1} holds if $\{b_l \mid 2 \leq l \leq n\}$ is linearly independent over $\mathbb{Z}$.  
\begin{hipotesis} \label{hip: 2}
    Every element in the set $\{b_l \mid 2 \leq l \leq n\}$ is equal to $1$. 
\end{hipotesis} 
\noindent Let $A \in \mathfrak{gl}(2n-1, \mathbb{R})$ as in \eqref{eq: matrixA} be given such that it satisfies either Hypothesis \ref{hip: 1} or \ref{hip: 2}. Clearly $\g_A$ is then completely solvable. Also, the $2$-form 
\begin{align} \label{eq: es casi la una de la mañana}
    \omega :=  e^1 \wedge e^{2n} + \sum_{i=2}^n e^i\wedge e^{i+n-1}
\end{align}
    \noindent turns out to be the symplectic form associated to an almost-K\"ahler structure on $\g_A$, and it is of central importance in this article; indeed, a non-trivial fact proved in this article is that if any Hypothesis holds then the associated Lefschetz operators $L_m^{(n)}:H^m(\g_A) \to H^{2n-m}(\g_A)$ in cohomology can be understood in terms of adjacency matrices of Kneser graphs of some parameters. 
We can now state the main result of the paper as follows: 
\begin{theorem} \label{thm: main result uno}
    Let $A \in \mathfrak{gl}(2n-1, \mathbb{R})$ as in \eqref{eq: matrixA} satisfy Hypothesis \ref{hip: 1}.  Then
\begin{enumerate}
    \item $G_A$ admits lattices $\Gamma$ for some choices of $b_2, \ldots, b_n$. 
    \item $H^*(\g_A) \cong H^*_{dR}(\Gamma \backslash G_A)$ can be computed explicitly, giving rise to natural ordered bases of both spaces. Their Betti numbers are given by
\begin{gather*}
    \dim H^{2k}(\g_A) = \dim H^{2k}_{dR}(\Gamma \backslash G_A) = \binom{n}{k}, \\
    \dim H^{2l+1}(\g_A) = \dim H^{2l+1}_{dR}(\Gamma \backslash G_A) = 2 \binom{n-1}{l}.
\end{gather*}
    \noindent for all $1\leq k \leq n$ and $0 \leq l \leq n-1$.   
    \item  For all $0 \leq m \leq n$, the matrix $M_m^{(n)}$ of the Lefschetz operator $L_m^{(n)}:H^m(\g_A) \to H^{2n-m}(\g_A)$ with respect to the symplectic form $\omega$ in \eqref{eq: es casi la una de la mañana} is the direct sum of adjacency matrices of Kneser graphs of some parameters for the choice of ordered bases of $H^m(\g_A)$ and $H^{2n-m}(\g_A)$ alluded to in (2). In particular, $L_m^{(n)}$ is invertible for all $0 \leq m \leq n$.
    \item $\g_A$ satisfies the hard-Lefschetz condition with respect to any symplectic form. Due to Hattori's theorem, the same holds for a solvmanifold $\Gamma\backslash G_A$ where $\Gamma$ is a lattice of $G_A$. 
\end{enumerate} 
\end{theorem}
\begin{theorem} \label{thm: main result dos}
    Let $A \in \mathfrak{gl}(2n-1, \mathbb{R})$ as in \eqref{eq: matrixA} satisfy Hypothesis \ref{hip: 2}.  Then
\begin{enumerate}
    \item $G_A$ admits lattices $\Gamma$. 
    \item $H^*(\g_A)\cong H^*_{dR}(\Gamma \backslash G_A)$ can be computed explicitly, giving rise to natural ordered bases of both spaces. Their Betti numbers are given by
\begin{gather*}
    \dim H^{2k}(\g_A) = \dim H^{2k}_{dR}(\Gamma \backslash G_A) = \binom{n-1}{k}^2 + \binom{n-1}{k-1}^2, \\
    \dim H^{2l+1}(\g_A) = \dim H^{2l+1}_{dR}(\Gamma \backslash G_A) = 2 \binom{n-1}{l}^2.
\end{gather*}
    \noindent for all $1\leq k \leq n$ and $0 \leq l \leq n-1$.   
    \item For all $0 \leq m \leq n$, the matrix $N_m^{(n)}$ of the Lefschetz operator $L_m^{(n)}:H^m(\g_A) \to H^{2n-m}(\g_A)$ with respect to the symplectic form $\omega$ in \eqref{eq: es casi la una de la mañana} is related to the matrices $M_m^{(n)}$ from (3) of Theorem \ref{thm: main result uno}, and is the direct sum of adjacency matrices of Kneser graphs of some parameters for the choice of ordered bases of $H^m(\g_A)$ and $H^{2n-m}(\g_A)$ alluded to in (2). In particular, $L_m^{(n)}$ is invertible for all $0 \leq m \leq n$.
    \item $\g_A$ satisfies the hard-Lefschetz condition with respect to any symplectic form. Due to Hattori's theorem, the same holds for a solvmanifold $\Gamma\backslash G_A$ where $\Gamma$ is a lattice of $G_A$. 
\end{enumerate} 
\end{theorem}
\indent Theorems \ref{thm: main result uno} and \ref{thm: main result dos} are generalizations of the main results in \cite{TT}. Item (4) of Theorems \ref{thm: main result uno} and \ref{thm: main result dos} are particular cases of the main result in \cite{K}. That result is established by a different approach than ours. We point out that our contributions are original in at least three different ways: all our constructions and computations are explicit, we establish a connection with graph theory, and we are able to find lattices for an appropriate choice of parameters.

\smallskip

\section{Preliminaries}

\subsection{The hard-Lefschetz condition} \label{seccion: the hard lefschetz condition}

\indent Let $(M,\omega)$ be a compact symplectic manifold of dimension $2n$, and denote by $\Omega(M)$ the algebra of differential forms in $M$. Informally, we say that $(M,\omega)$ satisfies the hard-Lefschetz condition, or simply that \textit{$(M, \omega)$ is hard-Lefschetz}, if it is possible to develop an analogue of Hodge theory for $(M,\omega)$. This loosely means that there is a notion of harmonicity of differential forms in $M$, depending on $\omega$ alone, such that every de Rham cohomology class in $M$ has a ``$\omega$-harmonic"\! representative. To clarify this, let us begin by noting that the $2$-form $\omega$ is closed and not exact by virtue of compactness and the Stokes theorem, and also that it induces a volume form $\vol_{\omega} = \frac{\omega^n}{n!}$ on $M$ due to its non-degeneracy, where $\omega^n \equiv \omega \wedge \cdots \wedge \omega$. The non-degeneracy condition is equivalent to the fact that the associated linear map 
\begin{align*}
    B:\mathfrak{X}(M) \to \mathfrak{X}(M)^* \equiv \Omega^1(M), \quad B(X)(Y) := \omega(X,Y)
\end{align*}
\noindent is bijective, and in particular that the skew-symmetric map 
\begin{align*}
    \omega^{-1}: \Omega(M) \times \Omega(M) \to \mathbb{R}, \quad \omega^{-1}(\theta, \eta) :=\eta( B^{-1}(\theta))
\end{align*}
\noindent is well defined. Therefore, $\omega$ induces a \textit{symplectic star operator} $\star_{\omega}:\Omega(M) \to \Omega(M)$ by the condition
\begin{align*}
    \alpha \wedge \star_{\omega} \beta = \omega^{-1}(\alpha,\beta) \vol_{\omega},
\end{align*}
\noindent from where it stems that it is a linear involution, and where $\omega^{-1}$ is the smooth bilinear form on $\Omega(M)$ given by the $C^\infty(M)$-linear extension of 
\begin{align*}
    \omega^{-1}( \theta_1 \wedge \cdots \wedge \theta_k, \eta_1 \wedge \cdots \wedge \eta_l ) = 
    \begin{cases}
        \det( \left[ \omega^{-1}(\theta_i, \eta_j) \right]_{1 \leq i,j \leq k} ), & k = l,\\
        0 & k \neq l,
    \end{cases}
\end{align*}
\noindent and $\theta_i$'s and $\eta_j$'s are $1$-forms. Note that $\star_{\omega}(\Omega^k(M)) \subseteq \Omega^{2n - k}(M)$ for all $0 \leq k \leq 2n$. The symplectic star operator $\star_{\omega}$ allows us to define a symplectic dual $d^c:\Omega(M) \to \Omega(M)$ of the exterior derivative $d$ by
\begin{align*}
   d^c\vert_{\Omega^k(M)} = (-1)^{k+1} \star_{\omega} \circ \, d \circ \star_{\omega},
\end{align*}
\noindent that certainly satisfies $d^c(\Omega^k(M)) \subseteq \Omega^{k-1}(M)$ for all $0 \leq k \leq 2n$ and $d^c \circ d^c = 0$. %If we denote by $d_k^c$ the restriction of $d^c$ to $\Omega^k(M)$, then it makes sense to consider the symplectic analogue of the de Rham cohomology:
%\begin{align*}
 %   H^k_{d^c}(M) := \frac{\operatorname{Ker} d_k^c}{\operatorname{Im} \; d_{k+1}^c}, \quad 0 \leq k \leq 2n; 
%\end{align*}
\noindent In particular, every $d^c$-exact form is $d^c$-closed.  A differential form $\alpha \in \Omega(M)$ is said to be \textit{$\omega$-harmonic} if it is $d$-closed and $d^c$-closed; that is, if $d \alpha = 0$ and $d^c \alpha = 0$.  

\indent We say that $(M, \omega)$ \textit{satisfies the $dd^c$-lemma} if every $d^c$-closed and $d$-exact form is $dd^c$-exact; in symbols, if any of the following equivalent identities hold:
\begin{align*}
    \operatorname{Ker} d^c \cap \mathrm{Im} \, d = \mathrm{Im} \, d d^c \quad \Longleftrightarrow \quad \operatorname{Ker} d \cap \mathrm{Im} \, d^c = \mathrm{Im} \, d d^c;
\end{align*}
\noindent the equivalence holds because $\star_\omega(\operatorname{Ker} d^c \cap \mathrm{Im} \; d) = \operatorname{Ker} d \cap \mathrm{Im} \; d^c$ and $dd^c = - d^c d$.

\indent We also have operators $L$, $\Lambda$, $H:\Omega(M) \to \Omega(M)$ given by
\begin{align*}
    L(\alpha) := \omega \wedge \alpha, \quad \Lambda(\alpha) := \star_{\omega} \circ L \circ \star_{\omega} \, \alpha, \quad H := [L, \Lambda].
\end{align*}
\noindent Note that, for all $0 \leq k \leq n$, we have
\begin{align*}
    L(\Omega^k(M)) \subseteq \Omega^{k+2}(M), \quad \Lambda(\Omega^k(M)) \subseteq \Omega^{k-2}(M), \quad H(\Omega^k(M)) \subseteq \Omega^k(M).
\end{align*}
\noindent These operators constitute a representation of the Lie algebra $\mathfrak{sl}(2, \mathbb{R})$ in $\Omega(M)$ identical to the one appearing in the well-known theory for compact K\"ahler manifolds. We mention in passing that $d^c = [d, \Lambda]$. The following operators obtained as repeated applications of the $L$-operator,
\begin{align*}
    L_k:\Omega^{n-k}(M) \to \Omega^{n+k}(M), \quad L_k(\alpha) = \omega^k \wedge \alpha,
\end{align*}
\noindent are of crucial importance to us, as they can be shown to induce operators in de Rham cohomology, which we call by the same name in a slight abuse of language as it is customary:
\begin{align}\label{eq: Lk}
    L_k: H_{dR}^{n-k}(M) \to H_{dR}^{n+k}(M), \quad L_k([\alpha]) = [\omega^k \wedge \alpha]. 
\end{align}
\noindent We say that $(M, \omega)$ is \textit{hard-Lefschetz} if the induced operators $L_k$ from \eqref{eq: Lk} are isomorphisms for all $0 \leq k \leq n$. We have the following result. 

\begin{theorem} \cite[Theorem 3.1]{TT} \label{thm: HL equivalences}
    The following conditions are equivalent on any compact symplectic manifold $(M, \omega)$ of dimension $2n$:
\begin{enumerate} [\rm (i)]
    \item $(M,\omega)$ is hard-Lefschetz. 
    \item Any de Rham cohomology class on $M$ has a $\omega$-harmonic representative.
    \item The $dd^c$-lemma holds. 
\end{enumerate}
\end{theorem} 
\noindent As indicated in \cite{TT}, Theorem \ref{thm: HL equivalences} is the culmination of works of several authors spanning many years of research. Other equivalent statements to those in Theorem \ref{thm: HL equivalences} can be found in the same reference. 

\indent The goal of this article is to find examples of \textit{solvmanifolds} that satisfy the hard-Lefschetz condition with respect to any invariant symplectic form. A review of solvmanifolds is in order. 

\smallskip

\subsection{Solvmanifolds} \label{seccion: solvmanifolds} Let $G$ be a connected real Lie group with Lie algebra $\g$. An almost complex structure $J$ on $G$ is left invariant if all left translations $L_g$, $g\in G$, are almost complex maps, that is, $(d L_g)_h J_h=J_{gh} (d L_g)_h$ for all $g,h\in G$. In this case $J$ is determined by the value at the identity of $G$. Thus, a left invariant complex structure on $G$ amounts to a  complex structure on its Lie algebra $\g$, that is, a real linear endomorphism $J$ of $\g$ satisfying $J^2 = -I$. A Riemannian metric $g$ on $G$ is called left invariant when all left translations are isometries. Such a metric $g$ is determined by its value $g_e=\pint$ at the identity $e$ of $G$, that is, $\pint$ is a positive definite inner product on $T_e G=\g$. 

An almost Hermitian structure $(J,g)$ on $G$ is called left invariant when both $J$ and $g$ are left invariant. Given such a left invariant almost Hermitian structure $(J,g)$ on $G$, let $J$ and $\pint$ denote the corresponding almost complex structure and almost Hermitian inner product on $\g$. We say that $(J, \pint)$ is an almost Hermitian structure on $\g$. If, moreover, the almost Hermitian structure $(J,g)$ is almost Kähler, i.e. the fundamental $2$-form $\omega=g(J\cdot,\cdot)$ is closed, then we also say that, at the Lie algebra level, $(J, \pint)$ is almost Kähler.

A \textit{solvmanifold} is a compact quotient $\Gamma\backslash G$, where $G$ is a simply connected solvable Lie group and $\Gamma$ is a discrete subgroup of $G$. Such a co-compact discrete subgroup $\Gamma$ is called a \textit{lattice} of $G$. When $G$ is nilpotent and $\Gamma\subset G$ is a lattice, the compact quotient $\Gamma\backslash G$ is known as a nilmanifold.

It follows that $\pi_1(\Gamma\backslash G)\cong \Gamma$ and  $\pi_n(\Gamma\backslash G)=0$ for $n>1$. Furthermore, solvmanifolds are determined up to diffeomorphism by their fundamental groups. In fact:

\begin{theorem}\cite{Mos}\label{thm: solv-isom}
If $\Gamma_1$ and $\Gamma_2$ are lattices in simply connected solvable Lie groups 
$G_1$ and $G_2$, respectively, and $\Gamma_1$ is isomorphic to $\Gamma_2$, then $\Gamma_1 \backslash G_1$ is diffeomorphic to $\Gamma_2 \backslash G_2$.
\end{theorem}

%The conclusion of the previous theorem can be strengthened when both solvable Lie groups G1G_1 and G2G_2 are completely solvable\footnote{A solvable Lie group GG is completely solvable if the adjoint operators \adx:\g→\g\ad_x:\g\to\g, with x∈\g=Lie(G)x\in \g=\operatorname{Lie}(G), have only real eigenvalues. In particular, nilpotent Lie groups are completely solvable.}. Indeed, this is the content of Saito's rigidity theorem:

%\begin{theorem}\cite{Sai}\label{Saito}
%Let G1G_1 and G2G_2 be simply connected completely solvable Lie groups and Γ1⊂G1,Γ2⊂G2\Gamma_1 \subset G_1, \, \Gamma_2\subset G_2 lattices. Then every isomorphism
%f:Γ1→Γ2f: \Gamma_1 \to \Gamma_2 extends uniquely to an isomorphism of Lie groups F:G1→G2F: G_1 \to G_2.
%\end{theorem}

A solvable Lie group $G$ is called completely solvable if the adjoint operators $\ad_x:\g\to\g$, with $x\in \g=\operatorname{Lie}(G)$, have only real eigenvalues. For solvmanifolds of completely solvable Lie groups the result in Theorem \ref{thm: solv-isom} can be strengthened, as the Saito's rigidity theorem shows:

\begin{theorem}\cite{Sai}\label{thm:Saito}
Let $G_1$ and $G_2$ be simply connected completely solvable Lie groups and $\Gamma_1 \subset G_1, \, \Gamma_2\subset G_2$ lattices. Then every isomorphism
$f: \Gamma_1 \to \Gamma_2$ 
extends uniquely to an isomorphism of Lie groups $F: G_1 \to G_2$.
\end{theorem}

Moreover, solvmanifolds of completely solvable Lie groups have a very useful property concerning their de Rham cohomology. Indeed, Hattori \cite{Hat} proved that the natural inclusion 
\[
\alt^* \g^* \hookrightarrow \Omega^*(\Gamma\backslash G),
\] 
with $G$ completely solvable, induces an isomorphism 
\begin{equation}\label{deRham}
H^*(\g) \cong H^*_{dR}(\Gamma\backslash G).
\end{equation}
That is, the de Rham cohomology of the solvmanifold can be computed in terms of left invariant forms. In particular, $H^*_{dR}(\Gamma\backslash G)$ does not depend on the lattice $\Gamma$.
The isomorphism \eqref{deRham} was proved earlier for nilmanifolds by Nomizu \cite{Nom}.  

In general, it is not easy to determine whether a given Lie group $G$ admits a lattice. A well known restriction is that if this is the case then $G$ must be unimodular (\cite{Mi}), i.e. the Haar measure on $G$ is left and right invariant, which is equivalent, when $G$ is connected, to  $\tr(\ad_x)=0$ for any $x$ in the Lie algebra $\g$ of $G$. In the nilpotent case there is a criterion to determine the existence of lattices, due to Malcev \cite{Mal}: a simply connected nilpotent Lie group has a lattice if and only if its Lie algebra admits a basis with respect to which the structure constants are rational.

We point out that left invariant geometric structures defined on $G$ induce corresponding geometric structures on $\Gamma\backslash G$, which are called invariant. For instance, a left invariant (almost) complex structure (respectively, Riemannian metric) on $G$ induces an (almost) complex structure (respectively, Riemannian metric) on $\Gamma\backslash G$ such that the canonical projection $G\to \Gamma\backslash G$ is a local biholomorphism (respectively, local isometry).

\subsection{Some graph theory}

\indent Perhaps surprisingly, the main results of the paper, Theorems \ref{thm: main result uno} and \ref{thm: main result dos}, are established using graph theory, so a very short review of this topic is in order. The results of this subsection are well-known. 

\begin{definition}
    A \textit{graph $X$}\footnote{We are excluding from the definition what in other places are called \textit{directed graphs}, \textit{graphs with loops}, and \textit{multigraphs}.} consists of a \textit{vertex set $V(X)$} and an \textit{edge set $E(X)$}, where an \textit{edge of $X$} is an unordered pair $\{x,y\}$ of distinct vertices of $X$. If $\{x,y\}$ is an edge then we say that $x$ and $y$ are \textit{adjoint vertices}. The \textit{adjacency matrix $A(X)$ of $X$} is a binary matrix with rows and columns indexed by the vertices of $X$ whose $(i,j)$-entry is $1$ if and only if $\{i,j\}$ is an edge of $X$, and zero otherwise. The graph $X$ is said to be \textit{invertible} if $A(X)$ is an invertible matrix.
\end{definition}
\begin{remark}
    Every binary symmetric matrix with zeros in the diagonal is the adjacency matrix of some graph. 
\end{remark}
\indent There is a family of graphs especially relevant in this article. 
\begin{definition}
    Fix $n$, $k \in \N$ with $n\geq 2$ and $0\leq k\leq n$. The \textit{Kneser graph of parameters $n$ and $k$}, denoted by $K(n,k)$, is the graph whose vertices are the subsets of $k$ elements of a set of $n$ elements, and where two vertices are adjacent if and only if the two corresponding subsets are disjoint.
\end{definition}
\indent Note that $K(n,k)$ has $\binom{n}{k}$ vertices and each vertex has \textit{degree} $\binom{n-k}{k}$; that is, every vertex in $K(n,k)$ is adjoint to precisely $\binom{n-k}{k}$ other vertices. The Kneser graph $K(n,1)$ is the \textit{complete graph with $n$ vertices}; that is, the graph in which every pair of vertices is an edge. We can interpret the graph $K(n,0)$ as being a collection of $n$ vertices with no edges. 
\begin{example} \label{ex: adjacency matrix of petersen graph}
    The adjacency matrix of the Kneser graph $K(5,2)$ is 
\begin{align*}
   A(K(5,2)) = 
   \left[
	\begin{array}{c c c c|c c c c c c}      
		0 & 0 & 0 & 0 & 0 & 0 & 0 & 1 & 1 & 1 \\
            0 & 0 & 0 & 0 & 0 & 1 & 1 & 0 & 0 & 1 \\
            0 & 0 & 0 & 0 & 1 & 0 & 1 & 0 & 1 & 0 \\
            0 & 0 & 0 & 0 & 1 & 1 & 0 & 1 & 0 & 0 \\
        \hline 
            0 & 0 & 1 & 1 & 0 & 0 & 0 & 0 & 0 & 1 \\
            0 & 1 & 0 & 1 & 0 & 0 & 0 & 0 & 1 & 0 \\
            0 & 1 & 1 & 0 & 0 & 0 & 0 & 1 & 0 & 0 \\
            1 & 0 & 0 & 1 & 0 & 0 & 1 & 0 & 0 & 0 \\
            1 & 0 & 1 & 0 & 0 & 1 & 0 & 0 & 0 & 0 \\
            1 & 1 & 0 & 0 & 1 & 0 & 0 & 0 & 0 & 0
	\end{array} \right]. 
\end{align*} 
    The ordering of the vertices is lexicographic; that is: 
\begin{align*}
    \{1,2\}, \quad \{1,3\}, \quad \{1,4\}, \quad\{1,5\}, \quad \{2,3\}, \quad \{2,4\}, \quad \{2,5\}, \quad \{3,4\}, \quad \{3,5\}, \quad \{4,5\}. 
\end{align*}
    \noindent Note that the upper-left all-zero block correspond to the vertices of the form $\{1, i\}$ for $2 \leq i \leq 5$. The block structure of the adjacency matrices of Kneser graphs is relevant in what follows (see Example \ref{ex: eme de petersen} below).
\end{example}
\indent Kneser graphs of all parameters are known to be invertible. Moreover, the spectrum of their adjacenty matrices is known, as the following result shows.

\begin{theorem} \label{thm: Kneser graphs are invertible} \cite[Theorem 9.4.3]{GR}. 
The eigenvalues of the Kneser graph $K(n,k)$ are the integers 
\begin{align*}
    \lambda_j = (-1)^j \binom{n-k-j}{k-j} \neq  0 \quad \text{for $0 \leq j \leq k$}. 
\end{align*}    
In particular, Kneser graphs (of all parameters) are invertible. 
\end{theorem} 
\indent The invertibility of Kneser graphs is crucial in what follows: it is closely related to the fact that the Lefschetz operators in the manifolds we study are linear isomorphisms. 

\smallskip 

\section{Diagonal almost Kähler almost abelian  Lie groups} \label{seccion: almost abelian lie groups}

\indent Let $A \in \mathfrak{gl}(d-1,\R)$ be given. Consider a real $d$-dimensional vector space $\g_A$ generated by $\{e_l \mid 1 \leq l \leq d\}$. We set
\begin{align*}
    \n := \Span_{\mathbb{R}} \{e_i \mid 1 \leq i \leq d-1 \}, \quad \n_0 := \Span_{\mathbb{R}} \{e_i \mid 2 \leq i \leq d-1 \},
\end{align*}
\noindent so as to have $\g_A := \R e_{d} \oplus \mathfrak{n}$. Define a bilinear skew-symmetric product 
$[\cdot, \cdot]: \g_A \times \g_A \to \g_A$ by setting 
\begin{align*}
    [\mathfrak{n}, \mathfrak{n}] = 0 \quad \text{and} \quad \text{$[e_{d}, v] = Av$ for all $v \in \n$}.  
\end{align*}
\noindent With these definitions, $\g_A$ turns out to be a 2-step solvable Lie algebra, with $\n$ a codimension one abelian ideal, and $\g_A$ is called an \textit{almost abelian} Lie algebra. Accordingly, the corresponding simply connected Lie group $G_A$ is a semidirect product $G_A=\R\ltimes_\phi \R^{d-1}$, where the action is given by $\phi(t)=\exp(t A)$. The Lie group $G_A$ is called almost abelian, as well. In fact, almost abelian Lie groups can be characterized as those Lie groups whose Lie algebra has a codimension one abelian ideal.

Note that $\g_A$ is completely solvable if and only if $A$ has real eigenvalues, and it is unimodular if and only if $\tr(A) = 0$. Moreover, $\n$ is the nilradical of $\g_A$ if and only if $A$ is not nilpotent (otherwise, $\g_A$ is nilpotent). It is known that $\g_A$ is isomorphic to $\g_{A'}$ if and only if $cA$ and $A'$ are conjugate for some scalar $c\neq 0$ (see \cite[Proposition 1]{Fr}). 

\smallskip

From now on we consider $d=2n$, so that $\g_A$ is $2n$-dimensional. If we consider on $\g_A$ the inner product $\langle \cdot, \cdot \rangle$ that makes $\{e_l \mid 1 \leq l \leq 2n\}$ an orthonormal basis and the orthogonal almost complex structure $J$ given by 
\begin{align} \label{eq: almost complex structure}
    %\omega = e^1 \wedge e^{2n} + \omega_0, \quad 
    J =
    \left[
	\begin{array}{c | c | c}      
		0 & 0 & -1 \\
            \hline
            0 & J_0 & 0 \\
            \hline
            1 & 0  & 0
	\end{array} \right] \quad \text{with} \quad 
    J_0 := 
    \left[
	\begin{array}{c | c}      
		0 & -\mathrm{Id} \\
            \hline
            \mathrm{Id} & 0
	\end{array} \right],
\end{align}
\noindent then the fundamental $2$-form $\omega$ is given by
\begin{equation}\label{eq: almost kahler structure}
    \omega = e^1 \wedge e^{2n} + \sum_{i=2}^n e^i\wedge e^{i+n-1},
\end{equation} 
%\noindent where $\omega_0=\sum_{i=2}^n e^i\wedge e^{i+n-1}\in \alt^2 \mathfrak{n}_0^*$ is a non-degenerate $2$-form in $\n_0^*$. 
\noindent When $d\omega=0$ the almost Hermitian structure $(J, \pint)$ is called almost K\"ahler.
Moreover, we have the following result.  

\begin{proposition} \label{prop: LW} \cite[Proposition 4.1, Remark 4.2]{LW}.
    Any almost Kähler almost abelian Lie algebra is equivalent to $(\g_A, J, g, \omega)$ with $J$ and $\omega$ as in \eqref{eq: almost complex structure} and \eqref{eq: almost kahler structure} and
\begin{align*}
    A = 
    \left[
	\begin{array}{c | c}      
		a & v^t \\
            \hline
            0 & A_0
	\end{array} \right] \quad \text{with} \quad
    A_0 := 
    \left[
	\begin{array}{c | c }      
		B & C \\
            \hline 
            D & - B^t 
	\end{array} \right],
\end{align*}
\noindent for some $a \geq 0$, $v \in \R^{2n-2}$, and $B$, $C$, $D \in \mathfrak{gl}(2n-2, \R)$ such that $C^t = C$ and $D^t = D$.     
\end{proposition}

\begin{remark}
    For the sake of completeness, we point out that using \cite[Lemma 6.1]{LRV} it is easily seen that the almost Hermitian structure in Proposition \ref{prop: LW} is K\"ahler, that is, $J$ is integrable, if and only if $v=0$ and $A_0J_0=J_0A_0$. This amounts to $B^t=-B$ and $D=-C$, i.e., $A_0\in\mathfrak{u}(n-1)$.  
\end{remark}

\indent In this article we deal with even-dimensional almost abelian Lie algebras with action given by a diagonal matrix $A \in \mathfrak{gl}(2n-1, \R)$. We are interested in the completely solvable and unimodular case, so we impose $a = 0$, $C = D = 0$ and $B = \diag(b_2, \ldots, b_n)$ for some $b_2, \ldots, b_n \in \R$. That is, we deal with $A$ of the form
\begin{equation}\label{eq: diagonal A}
    A = \diag(0,b_2, \ldots, b_n, -b_2, \ldots, -b_n).
\end{equation}
\noindent From now on we call $\g_A$ simply by $\g$.

We recall the following result regarding the existence and uniqueness of symplectic forms on diagonal almost abelian Lie algebras. 

\begin{proposition}\cite[Theorem 1.1]{CM}\label{prop:simpl-equiv}
Let $\g:=\R e_{2n}\ltimes_A \R^{2n-1}$ with $A=\operatorname{diag}(\lambda_1,\ldots, \lambda_{2n-1})$. Then:
\begin{enumerate}
\item There exists a symplectic form $\omega$ on $\g$ if and only if there exists a permutation $\pi$ of $\{1, \ldots, 2n-1\}$ such that $\lambda_{\pi(i)}+\lambda_{\pi(i+n)}=0$ for all $i$.
\item If there exists a symplectic form $\omega$ on $\g$, then it is unique up to
automorphism and scaling.
\end{enumerate}
\end{proposition}

\indent Recall that the standard cochain complex $(\alt^* \g^*, d)$ in any Lie algebra $\g$ with bracket $[\cdot, \cdot]$ has exterior derivative $d$ completely determined by the conditions
\begin{gather*}
    d\theta(x,y) = - \theta([x,y]), \quad \text{for all $\theta \in \g^*$ and any $x$, $y \in \g$}, \\
    d(\alpha \wedge \beta) = d \alpha \wedge \beta + (-1)^{k} \alpha \wedge d \beta, \quad \text{for all } \alpha\in\alt^k\g^*, \, \beta \in \alt^* \g^*.
\end{gather*}
\noindent Clearly, for $\g$ with $A$ as in equation \eqref{eq: diagonal A}, the exterior derivative $d:\g^*\to \alt^2\g^*$ is given by
\begin{equation}\label{eq:dif}
   d e^1=d e^{2n}=0, \quad d e^i=b_i e^{i} \wedge e^{2n} , \qquad d e^{i+n-1} = -b_i e^{i+n-1} \wedge e^{2n}, \quad 2\leq i\leq n.
\end{equation} 
This allows us to compute the exterior derivative of any $p$-form on~$\g$.
\begin{lemma} \label{lemma: exterior derivative formula}
    If $\beta \in \alt^\ast \mathfrak{n}_0^*$ is decomposable then there exist a scalar $c_{\beta}$, depending on $\beta$, such that
\begin{align*}
    d \beta = c_{\beta} \beta \wedge e^{2n}.
\end{align*}    
    \noindent More explicitly, if $\beta= e^{i_1}\wedge \cdots \wedge e^{i_p}$ for $2 \leq i_1 < \cdots < i_p \leq 2n-1 $ then
    \[ c_\beta= (-1)^{p+1} \left[ \sum_{l=1}^{\Lambda} b_{i_l} - \sum_{l = \Lambda +1}^p b_{i_l - n+1}  \right], \]
where $i_{\Lambda}$ is the largest index with value less or equal to $n$.
%\begin{align*}
%   d( e^{i_1} \wedge \cdots \wedge e^{i_{\Lambda}} \wedge e^{i_{\Lambda + 1}} \wedge \cdots \wedge e^{i_p} ) = (-1)^{p+1} \left[ \sum_{l=1}^{\Lambda} b_{i_l} - \sum_{l = \Lambda +1}^p b_{i_l - n+1}  \right] e^{i_1 \cdots i_{\Lambda} i_{\Lambda+1} \cdots i_p} \wedge e^{2n}.
%\end{align*}
%    \noindent Here, iΛi_{\Lambda} is the last index with value less or equal to nn. 
\end{lemma}

\begin{proof}
    The proof is by induction based on equation \eqref{eq:dif}. 
\end{proof}

\begin{corollary} \label{cor: del lema con la formula}
Let $\eta\in \alt^k\n_0^*$.
\begin{enumerate}
    \item If $\eta$ is of the form $\eta = \mu \wedge e^{2n}$ for some $(k-1)$-form $\mu$ then it is closed. 
    \item If $\eta$ is exact then $\eta = \mu \wedge e^{2n}$ for some $(k-1)$-form $\mu$.
    \item If $\eta$ satisfies $d\eta \neq 0$ then $\eta \wedge e^{2n}$ and $\delta \wedge \eta$ are exact. Here, $\delta := e^1 \wedge e^{2n}$.
\end{enumerate}
\end{corollary}
\begin{proof} \phantom{.} \\
\indent (1) and (2) follow immediately from Lemma \ref{lemma: exterior derivative formula}. 
   
(3) If $\eta$ is decomposable then, as $d \eta = c_{\eta} \eta \wedge e^{2n}$ because of Lemma \ref{lemma: exterior derivative formula}, the condition on $d \eta$ implies that $c_{\eta}$ is nonzero, from which it follows that $\eta \wedge e^{2n} = d \left( c_{\eta}^{-1} \eta \right)$. The case where $\eta$ is not decomposable is established by repeating the previous argument in every term of an appropriate expansion of $\eta$ into decomposable elements. Lastly, if $\eta \wedge e^{2n} = d \psi$ then $d (e^1 \wedge \psi) = \pm \delta \wedge \eta$.     
\end{proof}

Note that the kernel of the exterior derivative $d$ depends on how many of the expressions $\sum_{l=1}^{\Lambda} b_{i_l} - \sum_{l = \Lambda +1}^p b_{i_l - n+1}$ vanish. Our study will focus on two particular cases, in one of them few of these expressions vanish, while in the other many of them vanish.

\smallskip

\section{Computation of cohomology groups}
It is convenient to introduce some notation.  Let's fix for all the section $2 \leq i,j \leq n$ and, given $1 \leq k \leq n$, fix also $2 \leq i_1 < \cdots < i_k \leq n$ and $2 \leq j_1 < \cdots < j_k \leq n$. We begin by setting $\sigma(i) := i + n - 1$. This enables us to make sense of the expressions
\begin{align} \label{eq: thetas and gammas}
    \theta_{i|j} := e^i \wedge e^{\sigma(j)} \in \alt^2 \mathfrak{n}_0^*, \quad \gamma_i := \theta_{i \mid i} \in \alt^2 \mathfrak{n}_0^*. 
\end{align}
\noindent Note that every $\gamma_i$ is a closed $2$-form according to Lemma \ref{lemma: exterior derivative formula}.  We define in turn, for all $1 \leq k \leq n$,
\begin{gather*}
    \theta_{i_1 \cdots i_k \mid j_1 \cdots j_k} := \theta_{i_1 \mid j_1} \wedge \cdots \wedge \theta_{i_k \mid j_k} \in \alt^{2k} \mathfrak{n}_0^*, \\
    \gamma_{i_1 \cdots i_k} := \gamma_{i_1} \wedge \cdots \wedge \gamma_{i_k} \in \alt^{2k} \mathfrak{n}_0^*.
\end{gather*}
\noindent With this notation, we also have $\gamma_{i_1 \cdots i_k} = \theta_{i_1 \cdots i_k \mid i_1 \cdots i_k}$. It will become clear presently that the forms $\gamma_{i_1 \cdots i_k}$'s are fundamental in a sense that the forms $\theta_{i_1 \cdots i_k \mid j_1 \cdots j_k}$'s are not. To make sense of the following results, we recall the closed $2$-form $\delta$ from Corollary \ref{cor: del lema con la formula} given by
\begin{align*}
    \delta := e^1 \wedge e^{2n} \in \alt^* \mathfrak{g}^*.
\end{align*}     
\begin{lemma} \label{lemma: gammas cerradas no-exactas}
    For all $0 \leq k \leq n$, the following forms are closed and non-exact: 
\[\begin{array}{ll}
    \ri \, \gamma_{i_1 \cdots i_k} \in \alt^{2k} \mathfrak{n}_0^*, & \rii \, \delta \wedge \gamma_{i_1 \cdots i_k} \in \delta \wedge \alt^{2k} \mathfrak{n}_0^*, \\ &\\
    \riii \, e^1 \wedge \gamma_{i_1 \cdots i_k} \in \alt^{2k+1} \mathfrak{g}^*, & \riv \, \gamma_{i_1 \cdots i_k} \wedge e^{2n} \in \alt^{2k+1} \mathfrak{g}^*. 
\end{array}\]
    \noindent Moreover, nontrivial linear combinations of forms of each type are closed and non-exact as well. 
\end{lemma}

\begin{proof}
We know that the forms $e^1$, $e^{2n}$, and $\gamma_i$ are closed for $2\leq i \leq 2n-1$. This implies readily that forms of type $\ri$, $\rii$, $\riii$, and $\riv$, being wedge products of closed forms, are closed as well. Certainly, the same is true for linear combinations of those forms. Corollary \ref{cor: del lema con la formula} ensures that any linear combination of the forms of type $\ri$ or $\riii$ is not exact. Any nontrivial linear combination of forms of type $\riv$, say $\alpha := \sum c_{i_1 \cdots i_k}\gamma_{i_1 \cdots i_k}\wedge e^{2n}$, is non-exact. For if $\alpha = d \beta$ for some even form $\beta \in \alt^{2k} \mathfrak{n}_0^*$, say $\beta = \sum a_{i_1 \cdots i_{2k}} e^{i_1 \cdots i_{2k}}$ with $a_{i_1 \cdots i_{2k}} \neq 0$ and $e^{i_1 \cdots i_{2k}}$ nonclosed (and, in particular, not equal to a wedge product of $\gamma_i$'s), then from Lemma \ref{lemma: exterior derivative formula} we arrive at the absurd that $d \beta$ is not in the span of the $\gamma_{i_1 \cdots i_k} \wedge e^{2n}$'s. Any nontrivial linear combination of forms of type $\rii$, say $\alpha := \delta \wedge \sum c_{i_1 \cdots i_{k-1}}\gamma_{i_1 \cdots i_{k-1}}$, is also non-exact. To prove this, let us assume, by means of contradiction, that $\alpha = d \beta$ for some odd form $\beta \in \alt^{2k-1} \mathfrak{g}^*$. We can decompose $\beta$ as
\begin{align*}
    \beta = \beta_1 + e^1 \wedge \beta_2 + \beta_3 \wedge e^{2n} + \delta \wedge \beta_4, \qquad \beta_i\in \alt^*\n_0^*.
\end{align*}
    \noindent As both $\beta_3 \wedge e^{2n}$ and $\delta \wedge \beta_4 $ are closed by Corollary \ref{cor: del lema con la formula}, we may as well take $\beta_3$ and $\beta_4$ equal to zero without loss of generality. Likewise, as $d \beta_1 \in \alt^{2k-1} \mathfrak{n}_0^* \wedge e^{2n}$ (by Corollary \ref{cor: del lema con la formula} (2)) but $\iota_{e_1} (\alpha) \neq 0$, we must have that $\beta_1 = 0$: informally, $\beta_1$ is zero because ``$\alpha$ has an $e^1$"\! but ``$d \beta_1$ hasn't". Thus, we arrive at $\beta = e^1 \wedge \beta_2$. It follows from $\alpha = d (e^1 \wedge \beta_2)$ that $d\beta_2 = -\sum c_{i_1 \cdots i_{k-1}}\gamma_{i_1 \cdots i_{k-1}}\wedge e^{2n}$, and, again as a consequence of Lemma \ref{lemma: exterior derivative formula}, that $\beta_2$ is a linear combination of the forms $ \gamma_{i_1 \cdots i_{k-1}}$, and therefore that $\beta_2$ and $\beta$ are closed, from which it follows the absurd $\alpha = d \beta = 0$. 
\end{proof}
\indent Observe that we cannot assure which linear combinations of the forms $\theta_{i_1 \cdots i_k \mid j_1 \cdots j_k}$'s are closed without making any assumptions about the elements $b_1, \ldots, b_n$ of the matrix $A$ in \eqref{eq: diagonal A}. In this article we deal with two extreme cases, that can be described as to which one of the following hypotheses hold:

\begin{hypothesis} \label{hyp: 1}
    No finite nontrivial addition or subtraction of different elements in the set $\{b_l \mid 2 \leq l \leq n\}$ is zero. That is, if $\sum_{j=2}^n \varepsilon_jb_j=0$ for some $\varepsilon_j \in \{-1,0,1\}$ then $\varepsilon_j=0$ for all $j$.
\end{hypothesis}
    \indent Note that Hypothesis \ref{hyp: 1} holds if $\{b_l \mid 2 \leq l \leq n\}$ is linearly independent over $\mathbb{Z}$.  
\begin{hypothesis} \label{hyp: 2}
    Every element in the set $\{b_l \mid 2 \leq l \leq n\}$ is equal to $1$. 
\end{hypothesis} 
\indent The cases corresponding to Hypothesis \ref{hyp: 1} are in some sense the ones with the fewest closed non-exact $2$-forms, whereas the case corresponding to Hypothesis \ref{hyp: 2} is in some sense the one with the most closed non-exact $2$-forms. Corollaries \ref{cor: betti numbers of case I} and \ref{cor: betti numbers of case II} below offer a rigorous re-statement of this informal claim. In Section \ref{sec: lattices} we show that there are choices of $\{b_l \mid 2 \leq l \leq n\}$ that simultaneously satisfy Hypothesis \ref{hyp: 1} and the corresponding simply-connected Lie groups admit lattices; the same is true for the Lie group arising from Hypothesis \ref{hyp: 2}. 

\indent We introduce the following subspaces of $\alt^* \mathfrak{n}_0^*$: 
\begin{gather*}
    H^{2k}_{\mathrm{I}} := \Span_{\mathbb{R}} \{ \gamma_{i_1 \cdots i_k} \mid 2 \leq i_1 < \cdots < i_k \leq n \} \subseteq \alt^* \mathfrak{n}_0^*, \\
    H^{2k}_{\mathrm{II}} := \Span_{\mathbb{R}} \{ \theta_{i_1 \cdots i_k \mid j_1 \dots j_k} \mid 2 \leq i_1 < \cdots < i_k \leq n, \; 2 \leq j_1 < \cdots < j_k \leq n \} \subseteq \alt^* \mathfrak{n}_0^*. 
\end{gather*}
\noindent $H^{2k}_{\mathrm{I}}$ and $H^{2k}_{\mathrm{II}}$ will come in handy when describing cohomology groups of $\mathfrak{g}$ depending on which one of the Hypotheses above hold in $\mathfrak{g}$ (see Theorems \ref{thm: cohomology of case I} and \ref{thm: cohomology of case II} below).

\begin{lemma} \label{lemma: H_A cerradas}
    Assume that Hypothesis \ref{hyp: 1} holds.
\begin{enumerate}
    \item An even form $\alpha \in \alt^{2k} \mathfrak{n}_0^*$ is closed if and only if $\alpha \in H_{\mathrm{I}}^{2k}$. 
    \item An odd form $\alpha \in \alt^{2k+1} \mathfrak{n}_0^*$ is closed if and only if $\alpha=0$.
\end{enumerate}
\end{lemma} 
\begin{proof} \phantom{.} \\
    \indent (1) Let $\alpha = \sum_i a_i \beta_i$ be the expansion of $\alpha$ in linearly-independent decomposable elements of $\alt^{2k} \mathfrak{n}_0^*$, where we assume $a_i \neq 0$ for all $i$. Lemma \ref{lemma: exterior derivative formula} implies that $d \alpha = \sum_i a_i c_{\beta_i} \beta_i \wedge e^{2n}$. As the $\beta_i \wedge e^{2n}$ are linearly independent, then $\alpha$ is closed if and only if $c_{\beta_i} = 0$ for all $i$; according to Lemma \ref{lemma: exterior derivative formula} again, and using Hypothesis \ref{hyp: 1}, this happens if and only if each $\beta_i$ is a product of $\gamma_j$'s. \\
    \indent (2) Following the lines of the proof of (1), we may assume simply that $\alpha \in \alt^{2k+1}\n_0^*$ is a decomposable element. From $d\alpha=0$ and Lemma \ref{lemma: exterior derivative formula} we obtain that $c_\alpha=0$ and since the degree of $\alpha$ is odd, it follows again from Lemma \ref{lemma: exterior derivative formula}, using Hypothesis \ref{hyp: 1}, that this happens if and only $\alpha=0$. 
\end{proof}
\indent There is a refinement of Lemma \ref{lemma: gammas cerradas no-exactas} when Hypothesis \ref{hyp: 2} holds, as we state and prove right away. The validity of Hypothesis \ref{hyp: 1} does not lead to an obvious improvement of Lemma \ref{lemma: gammas cerradas no-exactas}. 
\begin{lemma} \label{lemma: thetas cerradas no-exactas}
    Assume that Hypothesis \ref{hyp: 2} holds. For all $0 \leq k \leq n$, the following forms are closed and non-exact: 
\[\begin{array}{ll}
    \ri \, \theta_{i_1 \cdots i_k \mid j_1 \cdots j_k} \in \alt^{2k} \mathfrak{n}_0^*, & \rii \, \delta \wedge \theta_{i_1 \cdots i_{k-1} \mid j_1 \cdots j_{k-1}} \in \delta \wedge \alt^{2k} \mathfrak{n}_0^*, \\ &\\
    \riii \, e^1 \wedge \theta_{i_1 \cdots i_k \mid j_1 \cdots j_k} \in \alt^{2k+1} \mathfrak{g}^*, & \riv \, \theta_{i_1 \cdots i_k \mid j_1 \cdots j_k} \wedge e^{2n} \in \alt^{2k+1} \mathfrak{g}^*. 
\end{array}\]
%\begin{gather*}
 %   \ri \, \theta_{i_1 \cdots i_k \mid j_1 \cdots j_k} \in \alt^{2k} \mathfrak{n}_0^*, \quad \rii \, \delta \wedge \theta_{i_1 \cdots i_{k-1} \mid j_1 \cdots j_{k-1}} \in \delta \wedge \alt^{2k} \mathfrak{n}_0^*, \\
  %  \riii \, e^1 \wedge \theta_{i_1 \cdots i_k \mid j_1 \cdots j_k} \in \alt^{2k+1} \mathfrak{n}_0^*, \quad \riv \, \theta_{i_1 \cdots i_k \mid j_1 \cdots j_k} \wedge e^{2n} \in \alt^{2k+1} \mathfrak{n}_0^*, 
%\end{gather*}   
    \noindent Moreover, nontrivial linear combinations of forms of each type are closed and non-exact as well. 
\end{lemma}
\begin{proof}
    Same as Lemma \ref{lemma: gammas cerradas no-exactas}, using Hypothesis \ref{hyp: 2} to establish that all the forms in the statement are closed. 
\end{proof}

\begin{lemma} \label{lemma: H_B cerradas}
    Assume that Hypothesis \ref{hyp: 2} holds.
\begin{enumerate}
    \item An even form $\alpha \in \alt^{2k} \mathfrak{n}_0^*$ is closed if and only if $\alpha \in H_{\mathrm{II}}^{2k}$. 
    \item An odd form $\alpha \in \alt^{2k+1} \mathfrak{n}_0^*$ is closed if and only if $\alpha=0$.
\end{enumerate}    
\end{lemma}
\begin{proof}
    Same as the proof of Lemma \ref{lemma: H_A cerradas}, using Hypothesis \ref{hyp: 2} instead of Hypothesis \ref{hyp: 1}. 
\end{proof} 
 \indent The cohomology groups of all orders in both cases of interest can be obtained as a direct consequence of this string of results. We introduce now some more notation to clear up the presentation. Denote by $\Gamma$ the top form in $\mathfrak{n}_0^*$ given by the product of all $\gamma_i$'s; that is, 
\begin{align*}
    \Gamma :=  \gamma_2 \wedge \cdots \wedge \gamma_n \in \alt^{2(n-1)} \mathfrak{n}_0^*.
\end{align*}
\noindent Similarly, denote by $\Gamma_{i_1 \cdots i_k}$ the product of all $\gamma_i$'s except for $\gamma_{i_1}, \ldots, \gamma_{i_k}$; that is, 
\begin{align*}
    \Gamma_{i_1 \cdots i_k} := \bigwedge_{l \in \{2, \ldots, n\} \setminus \{i_1, \ldots, i_k\}} \gamma_l \in \alt^{2(n-k-1)} \mathfrak{n}_0^*.
\end{align*}
\noindent By abuse of language we refer to any cohomology class by a suitably chosen representative of such class. In particular, from now on we take $H^{2k}_{\mathrm{I}}$ and $H^{2k}_{\mathrm{II}}$ to be subspaces of the cohomology group $H^{2k}$.

\begin{theorem} \label{thm: cohomology of case I}
    Assume Hypothesis \ref{hyp: 1} holds. The cohomology groups of $\g$ are given by
\begin{align*}  
    H^{2k} & = H^{2k}_{\mathrm{I}} \oplus \delta \wedge H^{2k-2}_{\mathrm{I}}, \quad 1\leq k\leq n, \\ 
    H^{2k+1} &= e^1 \wedge H^{2k}_{\mathrm{I}} \oplus H^{2k}_{\mathrm{I}} \wedge e^{2n}, \quad 0\leq k\leq n-1.
\end{align*}  
    \noindent In particular,
\begin{gather*}
    H^0 = \Span_{\mathbb{R}} \{1\}, \quad H^1 = \Span_{\mathbb{R}} \{e^1, e^{2n}\}, \\
    H^{2n-1} = \Span_{\mathbb{R}} \{e^1 \wedge \Gamma, \Gamma \wedge e^{2n}\}, \quad H^{2n} = \Span_{\mathbb{R}} \{\delta \wedge \Gamma\}. 
\end{gather*}   
\end{theorem}
\begin{proof}
    It follows from Lemma \ref{lemma: gammas cerradas no-exactas} that  
\begin{align*}
    H^{2k}_{\mathrm{I}} \subseteq H^{2k}, \quad \delta \wedge H^{2k-2}_{\mathrm{I}} \subseteq  H^{2k}, \quad e^1 \wedge H^{2k}_{\mathrm{I}} \subseteq H^{2k+1}, \quad H^{2k}_{\mathrm{I}} \wedge e^{2n} \subseteq H^{2k+1};
\end{align*}    
    \noindent as it is easily seen that $H^{2k}_{\mathrm{I}} \cap  \delta \wedge H^{2k-2}_{\mathrm{I}} = \{0\}$ and $ e^1 \wedge H^{2k}_{\mathrm{I}} \cap H^{2k}_{\mathrm{I}} \wedge e^{2n} = \{0\}$, it follows that
\begin{align*}
    H^{2k}_{\mathrm{I}} \oplus \delta \wedge H^{2k-2}_{\mathrm{I}} \subseteq H^{2k}, \quad e^1 \wedge H^{2k}_{\mathrm{I}} \oplus H^{2k}_{\mathrm{I}} \wedge e^{2n} \subseteq H^{2k+1}.
\end{align*}    
    \noindent Consider now a non-exact closed form $\beta$ of any degree. We decompose $\beta$ as
\begin{align*}
    \beta = \beta_1 + e^1 \wedge \beta_2 + \beta_3 \wedge e^{2n} + \delta \wedge \beta_4, \qquad \beta_i\in \alt^*\n_0^*;
\end{align*}
    \noindent notice that $\beta_1$ and $\beta_4$ share the same parity with $\beta$, whereas $\beta_2$ and $\beta_3$ have the opposite parity. As $e^1$ and $e^{2n}$ are closed $1$-forms, it follows from the formula given by Lemma \ref{lemma: exterior derivative formula} and the fact that $\beta$ is closed that
\begin{align*}
    0=d\beta= d\beta_1-e^1\wedge d\beta_2, 
\end{align*}    
    \noindent from where we get that $d\beta_1 = d \beta_2=0$, because $d\beta_1, \, d\beta_2\in \alt^*\n_0^*$. Corollary \ref{cor: del lema con la formula} implies that $d \beta_3 = d \beta_4 = 0$ as well: if $d \beta_3$ is nonzero then $\beta_3 \wedge e^{2n}$ is exact, and therefore zero in cohomology; the same goes for $\beta_4$. Thus,
\begin{itemize}
    \item if $\beta$ is an even form then, from Lemma \ref{lemma: H_A cerradas}, $\beta_1\in H^{2k}_{\mathrm{I}}$, $\beta_4 \in H^{2k-2}_{\mathrm{I}}$, and $\beta_2 = \beta_3 = 0$. 
    \item if $\beta$ is an odd form then, from Lemma \ref{lemma: H_A cerradas}, $\beta_2, \beta_3 \in H^{2k}_{\mathrm{I}}$, and $\beta_1 = \beta_4 = 0$. \qedhere 
\end{itemize}      
\end{proof}
\indent It is clear from the definition of $H^{2k}_{\mathrm{I}}$ that $\dim H^{2k}_{\mathrm{I}} = \binom{n-1}{k}$. As a consequence, we have the following result. 
\begin{corollary} \label{cor: betti numbers of case I}
    Assume Hypothesis \ref{hyp: 1} holds. The Betti numbers associated to $\mathfrak{g}$ are given by
\begin{gather*}
    b_{2k}(\g) = \binom{n}{k}, \quad  0\leq k \leq n, \\
    b_{2k+1}(\g) = 2 \binom{n-1}{k}, \quad  0\leq k \leq n-1.
\end{gather*}    
\end{corollary}
%%%
\begin{theorem} \label{thm: cohomology of case II}
Assume Hypothesis \ref{hyp: 2} holds. The cohomology groups of $\g$ are given by 
\begin{align*} 
    H^{2k} & = H^{2k}_{\mathrm{II}} \oplus \delta \wedge H^{2k-2}_{\mathrm{II}}, \quad 1\leq k\leq n, \\ 
    H^{2k+1} &= e^1 \wedge H^{2k}_{\mathrm{II}} \oplus H^{2k}_{\mathrm{II}} \wedge e^{2n}, \quad 0\leq k\leq n-1.
\end{align*}  
    \noindent In particular,
\begin{gather*}
    H^0 = \Span_{\mathbb{R}} \{1\}, \quad H^1 = \Span_{\mathbb{R}} \{e^1, e^{2n}\}, \\
    H^{2n-1} = \Span_{\mathbb{R}} \{e^1 \wedge \Gamma, \Gamma \wedge e^{2n}\}, \quad H^{2n} = \Span_{\mathbb{R}} \{\delta \wedge \Gamma\}. 
\end{gather*}    
\end{theorem}
\begin{proof}
    Same as in the proof of Theorem \ref{thm: cohomology of case I}, \textit{mutatis mutandis}. 
\end{proof}
\indent It is clear from the definition of $H^{2k}_{\mathrm{II}}$ that $\dim H^{2k}_{\mathrm{II}} = \binom{n-1}{k}^2$. As a consequence, we have the following result. 
\begin{corollary} \label{cor: betti numbers of case II}
    Assume Hypothesis \ref{hyp: 2} holds. The Betti numbers associated to $\mathfrak{g}$ are given by
\begin{gather*}
    b_{2k}(\g) = \binom{n-1}{k}^2 + \binom{n-1}{k-1}^2, \quad  1\leq k \leq n, \\
    b_{2k+1}(\g) = 2 \binom{n-1}{k}^2, \quad  0\leq k \leq n-1.
\end{gather*}   
\end{corollary}

\smallskip

\section{The hard-Lefschetz condition} 

In this section we verify that the hard-Lefschetz condition holds for any symplectic form on the almost abelian Lie algebras from the previous section, under the assumption of Hypothesis \ref{hyp: 1} or Hypothesis \ref{hyp: 2}. We divide the analysis in two cases, depending on which  hypothesis is considered.

\subsection{Case I}

\indent There are natural ordered bases for $H^{2k}$ and $H^{2k+1}$: assume that $1 \leq k,l \leq \floor{ \tfrac{n}{2} }$, and possibly that $k = 0$. Consider the sets\footnote{We adopt the convention that $B^{(n,0)} = C^{(n,1)} = \{1\}$.}
\begin{gather*}
    B^{(n,k)} := \{ \gamma_{i_1 \cdots i_k} \mid 2 \leq i_1 < \cdots < i_k \leq n\} \subseteq \alt^{2k} \mathfrak{n}_0^*, \\
    C^{(n,l)} := \{ \Gamma_{j_1 \cdots j_{l-1}} \mid 2 \leq j_1 < \cdots < j_{l-1} \leq n \} \subseteq \alt^{2(n-l)} \mathfrak{n}_0^*,
\end{gather*}
\noindent and define
\begin{gather*}
    B_1^{(n,k)} := \delta \wedge B^{(n,k-1)}, \quad B_2^{(n,k)} := B^{(n,k)}, \quad 
    B_3^{(n,k)} := e^1 \wedge B^{(n,k)}, \quad B_4^{(n,k)} := B^{(n,k)} \wedge e^{2n}, \\
    C_1^{(n,l)} := C^{(n,l)}, \quad C_2^{(n,l)} := \delta \wedge C^{(n,l+1)}, \quad
    C_3^{(n,l)} := e^1 \wedge C^{(n,l)}, \quad C_4^{(n,l)} := C^{(n,l)} \wedge e^{2n}.
\end{gather*}
\noindent The asymmetry in the definition of $B_1^{(n,k)}$ and $B_2^{(n,k)}$ with respect to $C_1^{(n,l)}$ and $C_2^{(n,l)}$ is intentional but ultimately inconvenient and will be dealt presently. For now, note that, according to Theorem \ref{thm: cohomology of case I},
\begin{multicols}{2}
    \begin{itemize}
        \item $B_1^{(n,k)} \cup B_2^{(n,k)}$ is a basis for $H^{2k}$.
        \item $B_3^{(n,k)} \cup B_4^{(n,k)}$ is a basis for $H^{2k+1}$.
        \item $C_1^{(n,l)} \cup C_2^{(n,l)}$ is a basis for $H^{2(n-l)}$.
        \item $C_3^{(n,l)} \cup C_4^{(n,l)}$ is a basis for $H^{2(n-l)+1}$.
    \end{itemize}
\end{multicols} 
\noindent Notice that the extreme cases $H^0$, $H^n$, and $H^{2n}$ are also taken into account by the previous description. We explicitly mention that the basis elements follow a lexicographic ordering in all cases. The following example should disambiguate what we mean.
\begin{example}
    Fix $n = 5$, and $k = l = 2$. The basis $B_1^{(5,2)} \cup B_2^{(5,2)}$ and $C_3^{(5,2)} \cup C_4^{(5,2)}$ of $H^{2k} = H^4$ and $H^{2(n-l)+1} = H^7$ respectively, written in the right order, are as follows:
\begin{gather*}
    B_1^{(5,2)} \cup B_2^{(5,2)} = \{ \delta \wedge \gamma_2, \delta \wedge \gamma_3, \delta \wedge \gamma_4, \delta \wedge \gamma_5 \} \cup \{ \gamma_{23}, \gamma_{24}, \gamma_{25}, \gamma_{34}, \gamma_{35}, \gamma_{45} \}, \\
    C_3^{(5,2)} \cup C_4^{(5,2)} = \{ e^1 \wedge \Gamma_2, e^1 \wedge \Gamma_3, e^1 \wedge \Gamma_4, e^1 \wedge \Gamma_5 \} \cup  \{ \Gamma_2 \wedge e^{2n}, \Gamma_3 \wedge e^{2n}, \Gamma_4 \wedge e^{2n}, \Gamma_5 \wedge e^{2n} \}.
\end{gather*}
    \noindent The ordering in the bases $B_3^{(5,2)} \cup B_4^{(5,2)}$ and $C_1^{(5,2)} \cup C_2^{(5,2)}$ are similar. According to our choice of notation, we have 
\begin{align*}
    C_3^{(5,2)} = \{ e^1 \wedge \gamma_{345}, e^1 \wedge \gamma_{245}, e^1 \wedge \gamma_{235}, e^1 \wedge \gamma_{234} \};
\end{align*}
\noindent these elements are written in the intended order. We acknowledge that it is better to adhere strictly to calling the elements of $C_1^{(n,l)}$, $C_2^{(n,l)}$, $C_3^{(n,l)}$, and $C_4^{(n,l)}$ using $\Gamma_{j_1, \ldots, j_{l-1}}$'s, and never using $\gamma_{i_1, \ldots, i_k}$'s. Not only this ensures that we always work with the intended ``forward lexicographic"\! ordering, but ultimately help us clarify the graph structure associated to the Lefschetz operators. 
\end{example} 
\indent The asymmetry in the definition of $B_1^{(n,k)}$ and $B_2^{(n,k)}$ with respect to $C_1^{(n,l)}$ and $C_2^{(n,l)}$ can be dealt with extending the range of indices to include $1$ and interpreting the $2$-form $\delta$ as a ``$\gamma_1$"\! of sorts, and redefining the elements $\Gamma_{j_1 \cdots j_s}$ accordingly. To avoid unnecessary confusions, call $\overline{\Gamma} := \delta \wedge \Gamma$, and
\begin{align} \label{eq: new notation}
    \overline{\gamma}_{i_1 i_2 \cdots i_k } :=
    \begin{cases}
        \delta \wedge \gamma_{i_2 \cdots i_l}  & i_1 = 1, \\
        \gamma_{i_1 i_2 \cdots i_l} & i_1 \geq 2,
    \end{cases} \quad 
    \overline{\Gamma}_{j_1 j_2 \cdots j_l} := 
    \begin{cases}
       \delta \wedge \Gamma_{j_1 j_2 \cdots j_l}  & j_1 \geq 2, \\
       \Gamma_{j_2 \cdots j_l} & j_1 = 1,
    \end{cases}
\end{align}
\noindent for all $1 \leq i_1 < i_2 < \cdots < i_k \leq n$ and $1 \leq j_1 < j_2 < \cdots < j_l \leq n$. Now define 
\begin{gather*}
    \overline{B}^{(n,k)} := \{ \overline{\gamma}_{i_1 \cdots i_k} \mid 1 \leq i_1 < \cdots < i_k \leq n\} \subseteq \alt^{2k} \mathfrak{g}^*, \\
    \overline{C}^{(n,l)} := \{ \overline{\Gamma}_{j_1 \cdots j_l} \mid 1 \leq j_1 < \cdots < j_l \leq n \} \subseteq \alt^{2(n-l)} \mathfrak{g}^*. 
\end{gather*}
\noindent Thus, we have
\begin{align*}
    \overline{B}^{(n,k)} = B_1^{(n,k)} \cup B_2^{(n,k)}, \quad \overline{C}^{(n,l)} = C_1^{(n,l)} \cup C_2^{(n,l)}.
\end{align*}
\noindent There is no reason to modify the description of the sets $B_3^{(n,k)}$, $B_2^{(n,k)}$, $C_3^{(n,k)}$, and $C_4^{(n,k)}$, so we leave them be. At the bottom of this anew description lie the bijections 
\begin{align*}
    B_1^{(n,k)} \longleftrightarrow C_1^{(n,k)}, \quad B_2^{(n,k)} \longleftrightarrow C_2^{(n,k)}, %\quad B_3^{(n,k)} \longleftrightarrow C_3^{(n,k+1)}, \quad B_4^{(n,k)} \longleftrightarrow C_4^{(n,k+1)}, 
\end{align*}
\noindent that now simply read
\begin{gather} \label{eq: new bijection I}  
    \overline{B}^{(n,k)} \ni \overline{\gamma}_{i_1 \cdots i_k} \quad \longleftrightarrow \quad \overline{\Gamma}_{i_1 \cdots i_k}  \in \overline{C}^{(n,k)} 
\end{gather}
\noindent for all $1 \leq i_1 < \dots < i_k \leq n$. Underlining the bijection in equation \eqref{eq: new bijection I} there is the fact that every choice of $m$ indices $r_1, \ldots, r_m$ in $\{1,2, \ldots, n\}$ corresponds to a complementary choice of indices $s_1, \ldots, s_{n-m}$, defined such that $\{r_1, \ldots, r_m\} \cup \{s_1, \ldots, s_{n-m}\} = \{1, 2, \ldots, n\}$. We stress that, while the bijection giving rise to \eqref{eq: new bijection I} corresponds to complementary choices in indices, two corresponding elements $\overline{\gamma}_{i_1 \cdots i_k}$ and $\overline{\Gamma}_{i_1 \cdots i_k}$ are represented by \underline{the same} choice $i_1, \ldots, i_k$, and will be done so in the rest of the article.  
\noindent We clarify all these ideas this with an example.
\begin{example} \label{ex: gamma 1}
    Fix $n = 5$ and $k = 2$. Consider the elements $\delta \wedge \gamma_4 \in B_1^{(5,2)}$ and $\gamma_{35} \in B_2^{(5,2)}$. We know that these bijectively correspond with the elements $\Gamma_4 = \gamma_{235} \in C_1^{(5,2)}$ and $\delta \wedge \Gamma_{35} = \delta \wedge \gamma_{24} \in C_2^{(5,2)}$, respectively. If we associate 
\begin{align*}
    \delta \wedge \gamma_4 \; \longleftrightarrow \; (1,4), \quad \gamma_{35} \; \longleftrightarrow \; (3,5), \quad \Gamma_4 \; \longleftrightarrow (2,3,5), \quad \delta \wedge \Gamma_{35} \; \longleftrightarrow \; (1,2,4),
\end{align*}
    \noindent then we can tell that $\delta \wedge \gamma_4$ and $\Gamma_4$ are associated to complementary choices of indices in the set $\{1, 2, 3, 4, 5\}$, and the same can be said about $\gamma_{35}$ and $\delta \wedge \Gamma_{35}$. If we call $\overline{\gamma}_{14} = \delta \wedge \gamma_4$, $\overline{\gamma}_{35} := \gamma_{35}$, $\overline{\Gamma}_{14} := \Gamma_4$, and $\overline{\Gamma}_{35} := \delta \wedge \Gamma_{35}$ then the correspondences are simply
\begin{align*}
    \overline{\gamma}_{14} \; \longleftrightarrow \; (1,4) \; \longleftrightarrow \; (2,3,5) \; &\longleftrightarrow \; \gamma_{235} = \overline{\Gamma}_{14}, \\
    \overline{\gamma}_{35} \; \longleftrightarrow \; (3,5) \; \longleftrightarrow \; (1,2,4) \; &\longleftrightarrow \; \delta \wedge \gamma_{24} = \overline{\Gamma}_{35}.  
\end{align*} 
\end{example}

\indent There are other bijections worth keeping in mind:
\begin{align*} 
    B_3^{(n,k)} \longleftrightarrow \overline{B}^{(n-1,k)} = B_1^{(n-1,k)} \cup B_2^{(n-1,k)}, \quad B_4^{(n,k)} \longleftrightarrow \overline{B}^{(n-1,k)} =  B_1^{(n-1,k)} \cup B_2^{(n-1,k)}, \\
    C_3^{(n,l)} \longleftrightarrow \overline{C}^{(n-1,l)} = C_1^{(n-1,l)} \cup C_2^{(n-1,l)}, \quad C_4^{(n,l)} \longleftrightarrow \overline{C}^{(n-1,l)} = C_1^{(n-1,l)} \cup C_2^{(n-1,l)}. 
\end{align*}
\noindent Each of these bijections is simply described by a shifting of every index by $1$; that is,
\begin{align} \label{eq: bijections I}
    B_3^{(n,k)} \ni e^1 \wedge \gamma_{i_1 i_2 \cdots i_k} \quad &\longleftrightarrow \quad \overline{\gamma}_{(i_1-1) (i_2-1) \cdots (i_k-1)} \in \overline{B}_3^{(n-1,k)}, \\
    C_4^{(n,l)} \ni \Gamma_{j_1 j_2 \cdots j_l} \wedge e^{2n} \quad &\longleftrightarrow \quad \overline{\Gamma}_{(j_1-1) (j_2-1) \cdots (j_l-1)} \in \overline{C}_4^{(n-1,l)}, \nonumber
\end{align}
\noindent and similarly for $B_4^{(n,k)}$ and $C_3^{(n,k)}$. We clarify this with an example. 
\begin{example}
    Fix $n = 6$ and $k = 2$.  Consider the elements $e^1 \wedge \gamma_{24}$ and $e^1 \wedge \gamma_{36}$ of $B_3^{(6,2)}$ in $H^5$, and $\Gamma_{25} \wedge e^{2n}$ and $\Gamma_{46} \wedge e^{2n}$ of $C_4^{(6,2)}$ in $H^7$. According to the bijections in equation \eqref{eq: bijections I}, all these elements correspond to
\begin{align*}
    B_3^{(6,2)} \ni e^1 \wedge \gamma_{24}    \quad &\longleftrightarrow \quad \overline{\gamma}_{13} = \delta \wedge \gamma_3 \in B_1^{(5,2)}, \\
    B_3^{(6,2)} \ni e^1 \wedge \gamma_{36}    \quad &\longleftrightarrow \quad \overline{\gamma}_{25} = \gamma_{25} \in B_2^{(5,2)}, \\
    C_4^{(6,2)} \ni \Gamma_{25} \wedge e^{2n} \quad &\longleftrightarrow \quad \overline{\Gamma}_{14} = \Gamma_4 \in C_1^{(5,2)}, \\
    C_4^{(6,2)} \ni \Gamma_{46} \wedge e^{2n} \quad &\longleftrightarrow \quad \overline{\Gamma}_{35} = \delta \wedge \Gamma_{35} \in C_2^{(5,2)}. 
\end{align*} 
\end{example}

\indent We now turn our attention to the hard-Lefschetz condition. We are interested in the $2$-form $\omega$ given by
\begin{align} \label{eq: argentina chile va 0 a 0}
    \omega = \delta + \sum_{l=2}^n \gamma_l = \sum_{l=1}^n \overline{\gamma}_l. 
\end{align}
\noindent By straightforward computations we arrive at the following result. 

\begin{lemma} \label{lemma: omega es no degenerada}
    For all $0 \leq k \leq n$ the $2$-form $\omega$ given by \eqref{eq: argentina chile va 0 a 0} satisfies
\begin{align*}
    \omega^k = k! \left[ \delta \wedge \sum_{2 \leq i_1 < \cdots < i_{k-1} \leq n} \gamma_{i_1 \cdots i_{k-1}} + \sum_{2 \leq j_1 < \cdots < j_k \leq n} \gamma_{j_1 \cdots j_k} \right] = k! \sum_{1 \leq i_1 < \cdots < i_k \leq n} \overline{\gamma}_{i_1 \cdots i_k}.
\end{align*}         
    \noindent In particular, for $k = n$, we have that $\omega^n = n! \; \delta \wedge \Gamma = n! \; \overline{\Gamma}$.
\end{lemma}
\indent Thus, $\omega$ is a symplectic form in $\mathfrak{g}$, being closed because of Lemma \ref{lemma: gammas cerradas no-exactas} and non-degenerate because of Lemma \ref{lemma: omega es no degenerada}. Note that $\omega$ is as in equation \eqref{eq: almost kahler structure}, and therefore is the symplectic form of an almost K\"ahler structure on $\g$ because of Proposition \ref{prop: LW}. According to Proposition \ref{prop:simpl-equiv}, $\omega$ is unique up to automorphism and scaling.  

\indent We now analyze the Lefschetz operators $L^{(n)}_m: H^m \to H^{2n-m}$ given by
\begin{align*}
    L^{(n)}_m(\cdot) := \frac{1}{(n-m)!} [\omega^{n - m} \wedge \cdot \;], \quad 0 \leq m \leq n.
\end{align*}
\noindent The dimension $n$ is stated explicitly because the structure of these operators depends on it. Notice that, similarly as before, we are indulging in the abuse of notation of calling ``$\omega$"\! to what should be called ``$[\omega]$". We have the following straightforward consequence of Lemma \ref{lemma: omega es no degenerada}; there we use the inner product $\langle \cdot, \cdot \rangle$ on $\alt^* \mathfrak{g}^*$ which is induced by the inner product on $\mathfrak{g}$ that makes $\{ e_l \mid 1 \leq l \leq n\}$ an orthonormal basis.  

\begin{proposition} \label{prop: properties of L I}
    For all $1 \leq k \leq \floor{ \tfrac{n}{2} }$, we have the following identities:
\begin{enumerate}
    \item 
\begin{align*}
    \langle L^{(n)}_{2k}( \delta \wedge \gamma_{i_1 \cdots i_{k-1}} ), \Gamma_{j_1 \cdots j_{k-1}}  \rangle = 0,  
\end{align*}
\begin{align*}
    \langle L^{(n)}_{2k}( \delta \wedge \gamma_{i_1 \cdots i_{k-1}} ), \delta \wedge \Gamma_{j_1 \cdots j_k}  \rangle = 
    \begin{cases}
        1, & \{i_1, \ldots, i_{k-1}\} \cap \{ j_1, \ldots, j_k \} = \emptyset, \\
        0, & \{i_1, \ldots, i_{k-1}\} \cap \{ j_1, \ldots, j_k \} \neq \emptyset. 
    \end{cases}
\end{align*}
\begin{align*}
    \langle L^{(n)}_{2k}( \gamma_{i_1 \cdots i_k} ), \Gamma_{j_1 \cdots j_{k-1}} \rangle = 
    \begin{cases}
        1, & \{i_1, \ldots, i_k\} \cap \{ j_1, \ldots, j_{k-1} \} = \emptyset, \\
        0, & \{i_1, \ldots, i_k\} \cap \{ j_1, \ldots, j_{k-1} \} \neq \emptyset.
    \end{cases}
\end{align*}
\begin{align*}
    \langle L^{(n)}_{2k}( \gamma_{i_1 \cdots i_k} ), \delta \wedge \Gamma_{j_1 \cdots j_k}  \rangle = 
    \begin{cases}
        1, & \{i_1, \ldots, i_k\} \cap \{ j_1, \ldots, j_k \} = \emptyset, \\
        0, & \{i_1, \ldots, i_k\} \cap \{ j_1, \ldots, j_k \} \neq\emptyset.
    \end{cases}
\end{align*}    
\noindent     
    \item 
\begin{align*}
    \langle L^{(n)}_{2k+1}( e^1 \wedge \gamma_{i_1 \cdots i_k} ), e^1 \wedge \Gamma_{j_1 \cdots j_k}  \rangle = 
    \begin{cases}
        1, & \{i_1, \ldots, i_k\} \cap \{ j_1, \ldots, j_k \} = \emptyset,  \\
        0, & \{i_1, \ldots, i_k\} \cap \{ j_1, \ldots, j_k \} \neq \emptyset.
    \end{cases}
\end{align*}
\begin{align*}
    \langle L^{(n)}_{2k+1}( e^1 \wedge \gamma_{i_1 \cdots i_k} ), \Gamma_{j_1 \cdots j_k} \wedge e^{2n} \rangle = 0,
\end{align*}    
\begin{align*}
    \langle L^{(n)}_{2k+1}( \gamma_{i_1 \cdots i_k} \wedge e^{2n} ), e^1 \wedge 
    \Gamma_{j_1 \cdots j_k} \rangle = 0,  
\end{align*}
\begin{align*}
    \langle L^{(n)}_{2k+1}( \gamma_{i_1 \cdots i_k} \wedge e^{2n}), \Gamma_{j_1 \cdots j_k} \wedge e^{2n} \rangle = 
    \begin{cases}
        1, & \{i_1, \ldots, i_k\} \cap \{ j_1, \ldots, j_k \} = \emptyset, \\
        0, & \{i_1, \ldots, i_k\} \cap \{ j_1, \ldots, j_k \} \neq \emptyset.
    \end{cases}
\end{align*}    
\end{enumerate}
\end{proposition} 

\indent The results in Proposition \ref{prop: properties of L I} become much cleaner if we use the formulas in equation \eqref{eq: new notation}, and even more so if we use multi-indices 
\begin{align*}
    I = (i_1, \ldots, i_k) &\text{ with $1 \leq i_1 < \cdots < i_k \leq n$}, \\
    J = (j_1, \ldots, j_k) &\text{ with $1 \leq j_1 < \cdots < j_k \leq n$}.  
\end{align*}
\noindent We say that the multi-indices $I$ and $J$ have \textit{length $k$}. We also express by $I \cap J = \emptyset$ the condition $\{i_1, \ldots, i_k\} \cap \{ j_1, \ldots, j_k \} = \emptyset$.
\begin{corollary} \label{cor: properties of L I}
    For all $1 \leq k \leq \floor{ \tfrac{n}{2} }$, $r$, $s \in \{1, 2n\}$, and multi-indices $I$ and $J$ of length k, the following identities hold. 
 \begin{gather*}
     \langle L_{2k}^{(n)}(\overline{\gamma}_I), \overline{\Gamma}_J \rangle = 
     \begin{cases}
         1, & I \cap J = \emptyset, \\ 
         0, & I \cap J \neq \emptyset,
     \end{cases} \\ 
     \langle L_{2k+1}^{(n)}(e^r \wedge \gamma_I), e^s \wedge \Gamma_J \rangle = 
     \begin{cases}
         1, & r = s \text{ and } I \cap J   =  \emptyset, \\
         0, & r \neq s \text{ or } I \cap J \neq \emptyset, \\ 
     \end{cases}
 \end{gather*}   
    \noindent In particular,
\begin{gather*}
    \langle L_{2k}^{(n)}(\overline{\gamma}_I), \overline{\Gamma}_J \rangle  = \langle L_{2k}^{(n)}(\overline{\gamma}_J), \overline{\Gamma}_I \rangle, \\
    \langle L_{2k+1}^{(n)}(e^1 \wedge \gamma_I), e^1 \wedge \Gamma_J \rangle = \langle L_{2k+1}^{(n)}(e^1 \wedge \gamma_J), e^1 \wedge \Gamma_I \rangle, \\
    \langle L_{2k+1}^{(n)}(\gamma_I \wedge e^{2n}), \Gamma_J \wedge e^{2n} \rangle = \langle L_{2k+1}^{(n)}(\gamma_J \wedge e^{2n}), \Gamma_I \wedge e^{2n} \rangle 
\end{gather*}    
    \noindent We are using the formulas in equation \eqref{eq: new notation}.
\end{corollary}
\indent It is clear from Corollary \ref{cor: properties of L I} that the structure of the Lefschetz operators in the odd case $m = 2k+1$ are, in some sense and barring the possible values of the indices, direct sum of two identical Lefschetz operators for the even case $m = 2k$. This was foreshadowed by the bijections in equation \eqref{eq: bijections I}. 

\indent Denote by $M_m^{(n)}$ the matrix of the operator $L_m^{(n)}$ with respect to the ordered basis of $H^m$ and $H^{2n-m}$. If we make the identifications leading up to equation \eqref{eq: new bijection I}, then we can also identify the entries of $M_m^{(n)}$ with a pair $(I,J)$ of ordered multi-indices of length $k$. 

\begin{remark}
    The fact that $B_1^{(n,k)} \cup B_2^{(n,k)}$ and $C_1^{(n,k)} \cup C_2^{(n,k)}$ are a basis for $H^{2k}$ and $H^{2(n-k)}$ respectively naturally gives the matrix $M^{(n)}_{2k}$ a block structure 
\begin{align*}
    M^{(n)}_{2k} = 
    \left[
	\begin{array}{c| c}      
		P & R \\
        \hline 
        S & Q 
	\end{array} \right],
\end{align*}
    \noindent where $P$ and $Q$ are square matrices of sizes $\binom{n-1}{k-1}$ and $\binom{n-1}{k}$, respectively. It is straightforward to see from Proposition \ref{prop: properties of L I} that $P \equiv 0$ and that $Q$ is antidiagonal, although the latter is somewhat tedious to prove. We won't use these properties of $M^{(n)}_{2k}$ in the sequel. What else can be said about the matrices $M^{(n)}_{2k}$ is already in Lemma \ref{lemma: properties of L I}, and so do every interesting feature of the matrices $M^{(n)}_{2k+1}$. 
\end{remark}

\begin{example} \label{ex: eme de petersen}
    For $n = 5$ and $k = 2$ we have 
\begin{align*}
    M^{(5)}_4 = \left[
	\begin{array}{c c c c|c c c c c c}      
		0 & 0 & 0 & 0 & 0 & 0 & 0 & 1 & 1 & 1 \\
        0 & 0 & 0 & 0 & 0 & 1 & 1 & 0 & 0 & 1 \\
        0 & 0 & 0 & 0 & 1 & 0 & 1 & 0 & 1 & 0 \\
        0 & 0 & 0 & 0 & 1 & 1 & 0 & 1 & 0 & 0 \\
        \hline 
        0 & 0 & 1 & 1 & 0 & 0 & 0 & 0 & 0 & 1 \\
        0 & 1 & 0 & 1 & 0 & 0 & 0 & 0 & 1 & 0 \\
        0 & 1 & 1 & 0 & 0 & 0 & 0 & 1 & 0 & 0 \\
        1 & 0 & 0 & 1 & 0 & 0 & 1 & 0 & 0 & 0 \\
        1 & 0 & 1 & 0 & 0 & 1 & 0 & 0 & 0 & 0 \\
        1 & 1 & 0 & 0 & 1 & 0 & 0 & 0 & 0 & 0
	\end{array} \right]. 
\end{align*}    
    \noindent The block structure of $M_4^{(5)}$ is shown explicitly. Examples for larger $n$ and $k$ aren't particularly illuminating because of the  enormous size of $M_4^{(5)}$, so we omit them. Notice nonetheless that $M_4^{(5)}$ is the adjacency matrix of the Kneser graph $K(5,2)$ as per Example \ref{ex: adjacency matrix of petersen graph}, a fact that is decidedly not a coincidence.   
\end{example}

\begin{lemma} \label{lemma: properties of L I} 
    The following statements hold for all $1 \leq k \leq \floor{ \tfrac{n}{2} }$. 
\begin{enumerate} 
    \item $M_{2k}^{(n)}$ is a binary symmetric matrix with zeros in the diagonal. The $(I,J)$-entry of $M_{2k}^{(n)}$ is nonzero if and only if $I \cap J = \emptyset$. 
    \item $M_{2k+1}^{(n)}$ is of the form $\left(
	\begin{array}{c|c}      
		M_{2k}^{(n-1)} & 0 \\
        \hline 
        0 & M_{2k}^{(n-1)} \\
	\end{array} \right)$. In particular, $M_{2k+1}^{(n)}$ is a binary symmetric matrix with zeros in the diagonal.   
\end{enumerate}    
\end{lemma}
\begin{proof} \phantom{.} \\
    \indent (1) It is clear from Corollary \ref{cor: properties of L I} that $M_{2k}^{(n)}$ is a binary symmetric matrix, and that the $(I,J)$-entry of $M_{2k}^{(n)}$ is nonzero if and only if $I \cap J = \emptyset$. An element in the diagonal of $M_{2k}^{(n)}$ is either
\begin{align*}
    \langle L^{(n)}_{2k}( \delta \wedge \gamma_{i_1 \cdots i_{k-1}} ), \Gamma_{i_1 \cdots i_{k-1}}  \rangle \quad \text{or} \quad \langle L^{(n)}_{2k}( \gamma_{i_1 \cdots i_k} ), \delta \wedge \Gamma_{i_1 \cdots i_k} \rangle  
\end{align*}    
    \noindent for some $2 \leq i_2 < \cdots < i_k \leq n$, and all those elements are zero because of Proposition \ref{prop: properties of L I}. 
    \indent (2) The block-diagonal form of $M^{(n)}_{2k+1}$ follows from the identities 
\begin{align*}
    \langle L^{(n)}_{2k+1}( e^1 \wedge \gamma_{i_1 \cdots i_k} ), \Gamma_{j_1 \cdots j_k} \wedge e^{2n} \rangle = 0, \quad 
    \langle L^{(n)}_{2k+1}( \gamma_{i_1 \cdots i_k} \wedge e^{2n} ), e^1 \wedge 
    \Gamma_{j_1 \cdots j_k} \rangle = 0,  
\end{align*}
    \noindent appearing in Proposition \ref{prop: properties of L I}. The fact that each block is precisely $M^{(n-1)}_{2k+1}$ follows from the bijections in equation \eqref{eq: bijections I}. \qedhere    
\end{proof}

\indent Lemma \ref{lemma: properties of L I} allows us to interpret the matrices $M^{(n)}_m$ as the adjacency matrices of some graphs. Indeed, as the $(I,J)$-entry of $M^{(n)}_{2k}$ is nonzero if and only if $I \cap J = \emptyset$, it is the adjacency matrix of the Kneser graph $K(n,k)$ of parameters $n$ and $k$. In particular, as per Theorem \ref{thm: Kneser graphs are invertible}, $M^{(n)}_{2k}$ is invertible, so that the Lefschetz operator $L^{(n)}_{2k}$ is a linear isomorphism. By the same vein, $M^{(n)}_{2k+1}$ is the adjacency matrix of a graph with two connected components, where each component is a Kneser $K(n-1,k)$ graph of parameters $n-1$ and $k$.
\begin{remark} \label{obs: el caso con k = 0}
    The case with $k = 0$ and arbitrary $n$ can be dealt with easily: Given that $L_0(\cdot) = \frac{1}{n!} \omega^n \wedge \cdot = \delta \wedge \Gamma \wedge \cdot$, $M_0^n$ is simply the identity matrix. 
\end{remark}
\indent We have established the main result of the subsection. 

\begin{theorem} \label{thm: main result I}
    The Lie algebras $\g = \R e_{2n}\ltimes_A \R^{2n-1}$, with $A$ given by \eqref{eq: diagonal A} where Hypothesis \ref{hyp: 1} holds, satisfy the hard-Lefschetz condition with respect to any symplectic form $\omega$ on $\g$.
\end{theorem}
\begin{proof}
The previous paragraph allows us to establish the result for $\omega = \delta + \sum_{l=2}^n \gamma_l$, as the adjacency matrix of the Kneser graphs $K(n,k)$ and $K(n-1,k)$ are invertible because of Theorem \ref{thm: Kneser graphs are invertible}, and treating the case $k = 0$ separately. The stated result follows from this because of Proposition \ref{prop:simpl-equiv}(ii).
\end{proof}

\smallskip

\subsection{Case II}
As in Case I, there are natural ordered bases for $H^{2k}$ and $H^{2k+1}$, albeit slightly more complicated. Recall the definition of the elements $\theta_{i|j}$ from equation \eqref{eq: thetas and gammas}. For clarity, in this subsection we will \underline{never} write $\theta_{i|j}$ if $i = j$; that is, no $\theta_{i|j}$ is a $\gamma_l$. Denote the set of ordered multi-indices of length $m$ of the set $\{2, \ldots, n\}$ by $\mathcal{I}_m^{(n)}$, and the set of ordered multi-indices of length $m$ of the set $\{1, 2, \ldots, n\}$ by $\overline{\mathcal{I}}_m^{(n)}$.

\indent We describe the bases of $H^{2k}$ and $H^{2k+1}$ as follows. Assume that $1 \leq k,l \leq \floor{ \tfrac{n}{2} }$, and possibly that $k = 0$. Fix $0 \leq p  \leq k$ and $0 \leq q  \leq l$. Pick  
\begin{gather*}
    I = (i_1, \ldots, i_{k-p}) \in \mathcal{I}_{k-p}^{(n)}, \quad J = (j_1, \ldots, j_{l-q}) \in \mathcal{I}_{l-q}^{(n)}, \\ 
    R_p = (r_1, \ldots, r_p) \in \mathcal{I}_p^{(n)}, \quad S_p = (s_1, \ldots, s_p) \in \mathcal{I}_p^{(n)}, \\
    U_q = (u_1, \ldots, u_q) \in \mathcal{I}_q^{(n)}, \quad V_q = (v_1, \ldots, v_q) \in \mathcal{I}_q^{(n)},
\end{gather*}
\noindent such that the following conditions hold:
\begin{gather} \label{eq: the conditions}
    R_p \cap S_p = \emptyset, \quad \emptyset = U_q \cap V_q, \\
    I \cap (R_p \cup S_p) = \emptyset, \quad \emptyset = (U_q \cup V_q) \cap J. \notag 
\end{gather}
\noindent Note that the slightly different notation is intentional. Recall equation \eqref{eq: new notation}. We define
\begin{gather*}
    B^{(n,k)}_{(R_p,S_p)} := \{ \gamma_I \wedge \theta_{R_p | S_p} \mid I \in \mathcal{I}_{k-p}^{(n)} \}, \quad C^{(n,l)}_{(U_q,V_q)} := \{ \Gamma_J \wedge \theta_{U_q | V_q} \mid J \in \mathcal{I}_{l-q}^n \}, \\ 
    \overline{B}^{(n,k)}_{(R_p,S_p)} := \{ \overline{\gamma}_I \wedge \theta_{R_p | S_p} \mid I \in \overline{\mathcal{I}}_{k-p}^{(n)} \}, \quad \overline{C}^{(n,l)}_{(U_q,V_q)} := \{ \overline{\Gamma}_J \wedge \theta_{U_q | V_q} \mid J \in \overline{\mathcal{I}}_{l-q}^n \}.
\end{gather*}
\noindent We also define 
\begin{gather*}
    \overline{B}^{(n,k)} := \bigcup_{p = 0}^k \bigcup_{(R_p,S_p)}\overline{B}^{(n,k)}_{(R_p,S_p)}, \quad \overline{C}^{(n,l)} := \bigcup_{q = 0}^l \bigcup_{(U_q, V_q)} \overline{C}^{(n,l)}_{(U_q,V_q)}, \\
    B_3^{(n,k)} := \bigcup_{p = 0}^k \bigcup_{(R_p,S_p)} e^1 \wedge B^{(n,k)}_{(R_p,S_p)}, \quad C_3^{(n,l)} := \bigcup_{q = 0}^l \bigcup_{(U_q,V_q)} e^1 \wedge C^{(n,l)}_{(U_q,V_q)}\\
    B_4^{(n,k)} := \bigcup_{p = 0}^k \bigcup_{(R_p,S_p)} B^{(n,k)}_{(R_p,S_p)} \wedge e^{2n}, \quad C_4^{(n,l)} := \bigcup_{q = 0}^l \bigcup_{(U_q,V_q)} C^{(n,l)}_{(U_q,V_q)} \wedge e^{2n}.
\end{gather*}
\noindent The subscripts $3$ and $4$ stand only to reminisce the corresponding definitions of Case I. Note that, according to Theorem \ref{thm: cohomology of case II},
\begin{multicols}{2}
    \begin{itemize}
        \item $\overline{B}^{(n,k)}$ is a basis for $H^{2k}$.
        \item $B_3^{(n,k)} \cup B_4^{(n,k)}$ is a basis for $H^{2k+1}$.
        \item $\overline{C}^{(n,l)}$ is a basis for $H^{2(n-l)}$.
        \item $C_3^{(n,l)} \cup C_4^{(n,l)}$ is a basis for $H^{2(n-l)+1}$.
    \end{itemize}
\end{multicols} 
\noindent  As in Case I, notice that the extreme cases $H^0$, $H^n$, and $H^{2n}$ are also taken into account by the previous description. The order in these bases is always lexicographic\footnote{The preferred order is for increasing $p$, following first the lexicographic order of the indices $I$'s, then the $R$'s, and lastly the $S$'s. }, although what really matters is that only the product of $\gamma_l$'s and $\Gamma_j$'s (or $\overline{\gamma}_l$'s and $\overline{\Gamma}_j$'s) are ordered in that way; we'll see presently that the ordering concerning the indices corresponding to the $\theta_{i|j}$'s is immaterial.  
\begin{remark}
    According to our conventions, the element $\Gamma_{i_1 \cdots i_{k-p}} \wedge \theta_{r_1 \cdots r_p | s_1 \cdots s_p}$ is a wedge product of $\gamma_l$'s whose corresponding indices are \textit{not} one of $i_1, \ldots, i_{k-p}$, and a wedge product of $\theta_{i|j}$'s whose corresponding indices \textit{are} one of the pairs $(r_k, s_k)$. 
\end{remark} 

\indent Similarly as in case I, note that  for every $0 \leq p \leq k$ and every pair of multi-indices $(R_p,S_p)$ there are bijections 
\begin{align} \label{eq: new bijection II}  
    \overline{B}_{(R_p,S_p)}^{(n,k)} \longleftrightarrow \overline{C}_{(R_p,S_p)}^{(n,k)},
\end{align}
\noindent and 
\begin{gather} \label{eq: bijections II} 
    e^1 \wedge B_{(R_p,S_p)}^{(n,k)} \longleftrightarrow \overline{B}_{(R_p,S_p)}^{(n-1,k)}, \quad B_{(R_p,S_p)}^{(n,k)} \wedge e^{2n} \longleftrightarrow \overline{B}_{(R_p,S_p)}^{(n-1,k)}, \\
    e^1 \wedge C_{(R_p,S_p)}^{(n,k)} \longleftrightarrow \overline{C}_{(R_p,S_p)}^{(n-1,k)}, \quad C_{(R_p,S_p)}^{(n,k)} \wedge e^{2n} \longleftrightarrow \overline{C}_{(R_p,S_p)}^{(n-1,k)}, \nonumber
\end{gather}
\noindent defined in an analogous way to the bijections given in equations \eqref{eq: new bijection I} and \eqref{eq: bijections I}.   

\indent We now turn our attention to the hard-Lefschetz condition. We use the same symplectic form $\omega$ as in Case I,
\begin{align*}
    \omega = \delta + \sum_{l=2}^n \gamma_l = \sum_{l=1}^n \overline{\gamma}_l. 
\end{align*}
\noindent In particular, Lemma \ref{lemma: omega es no degenerada} still holds.  What's really convenient about this symplectic form is that 
\begin{align} \label{eq: omega preserva theta}
    \omega^m \wedge B^{(n,k)}_{(R_p,S_p)} \subseteq \overline{B}^{(n,k+m)}_{(R_p,S_p)} \quad \text{for all $m \in \mathbb{N}$ and $0 \leq p \leq k$},
\end{align}
\noindent a condition that can be described heuristically as ``it doesn't change $p$"\! or ``it preserves the amount of $\theta_{i|j}$'s". This feature indicates the direct-sum structure that characterize Case II, and ultimately motivates the ugly description of the bases of $H^{2k}$ and $H^{2k+1}$.

\indent We now analyze the Lefschetz operators $L^{(n)}_m: H^m \to H^{2n-m}$ given by
\begin{align*}
    L^{(n)}_m(\cdot) := \frac{1}{(n-m)!} [\omega^{n - m} \wedge \cdot \;], \quad 0 \leq m \leq n.
\end{align*}
\noindent Notice that, similarly as before, we are indulging in the abuse of notation of calling ``$\omega$"\! to what should be called ``$[\omega]$". We have the  following corollary of Lemma \ref{lemma: omega es no degenerada} and equation \eqref{eq: omega preserva theta}; we'll skip the result analogous to Proposition \ref{prop: properties of L I}, and jump right away to the cleaner version of it corresponding to Corollary \ref{cor: properties of L I}.

\begin{corollary} \label{cor: properties of L II}
    For all $1 \leq k \leq \floor{ \tfrac{n}{2} }$, $r$, $s \in \{1, 2n\}$, and multi-indices $I$, $J$, $R_p$, $S_p$, $U_q$, $V_q$ that satisfy the conditions given by equation \eqref{eq: the conditions}, the following identities hold.  
 \begin{gather*}
     \langle L_{2k}^{(n)}(\overline{\gamma}_I \wedge \theta_{R_p | S_p}), \overline{\Gamma}_J \wedge \theta_{U_q | V_q} \rangle = 
     \begin{cases}
         1, & p = q, \; I \cap J = \emptyset, \; R_p = U_q, \; S_p = V_q, \\
         0, & \text{otherwise}. 
     \end{cases} \\ 
     \langle L_{2k+1}^{(n)}(e^r \wedge \gamma_I \wedge \theta_{R_p | S_p} ), e^s \wedge \Gamma_J  \wedge \theta_{U_q | V_q}  \rangle = 
     \begin{cases}
         1, & r = s, \; p = q, \; I \cap J = \emptyset, \; R_p = U_q, \; S_p = V_q, \\
         0, & \text{otherwise}.
     \end{cases}
 \end{gather*}   
    \noindent We are using the formulas in equation \eqref{eq: new notation}.
\end{corollary}

\indent Denote by $N_m^{(n)}$ the matrix of the operator $L_m^{(n)}$ with respect to the ordered basis of $H^m$ and $H^{2n-m}$. Similar considerations to those in Case I apply in this case as well. The block structure of the matrices $N_m^{(n)}$ can be obtained in a straightforward manner from Corollary \ref{cor: properties of L II} much in the same way as was done in Case I. 

\begin{lemma} \label{lemma: properties of L II}
The following statements hold for all $0 \leq k \leq \floor{ \tfrac{n}{2} }$. 
\begin{enumerate} 
    \item $N_{2k}^{(n)}$ is a block-diagonal matrix of the form 
\begin{align*}
    N_{2k}^{(n)} = \bigoplus_{p=0}^k \bigoplus_{(R_p,S_p)} M^{(n-2p)}_{2(k-p)},
\end{align*}
    \noindent where $M^{(n-2p)}_{2(k-p)}$ is a matrix as in Lemma \ref{lemma: properties of L I}. In particular, $N_{2k}^{(n)}$ is similar to a block diagonal matrix, with one block being binary symmetric having zeros in the diagonal and the other being an identity matrix of size $k \times k$. 
    \item $N_{2k+1}^{(n)}$ is a block-diagonal matrix of the form 
\begin{align*}
    N_{2k+1}^{(n)} = \bigoplus_{p=0}^k \bigoplus_{(R_p,S_p)} M^{(n-2p)}_{2(k-p)+1},
\end{align*}
    \noindent where $M^{(n-2p)}_{2(k-p)}$ is a matrix as in Lemma \ref{lemma: properties of L I}. In particular, $N_{2k+1}^{(n)}$ is similar to a block diagonal matrix, with one block being binary symmetric having zeros in the diagonal and the other being an identity matrix of size $2k \times 2k$. 
\end{enumerate}    
\end{lemma}
\begin{proof}
    The direct-sum structure follows from equation \eqref{eq: omega preserva theta}. Every summand is as stated because of Corollary \ref{cor: properties of L II}, which essentially ensures that, for every fixed $\theta_{U_p \mid V_p}$, $N_m^{(n)}$ behaves exactly as in Case I but with the $2p$ indices $\{r_1, \ldots, r_p, s_1, \ldots, s_p\}$ removed. The existence of a identity-block of the stated size follows from the fact that for every $k = p$ the Lefschetz operator is the identity operator, as was pointed out in Remark \ref{obs: el caso con k = 0}. 
\end{proof} 

\indent As in Case I, Corollary \ref{cor: properties of L II} and Lemma \ref{lemma: properties of L II} allow us to interpret the matrices $N^{(n)}_m$ as direct sum of adjacency matrices of several Kneser graphs of parameters $n$ and $k - p$, with $0 \leq p \leq k$, and a identity matrix. Every one of these blocks in invertible, partly because of Theorem \ref{thm: Kneser graphs are invertible}, and so we have established the main result of the subsection. 

\begin{theorem} \label{thm: main result II}
    The Lie algebras $\g = \R e_{2n}\ltimes_A \R^{2n-1}$, with $A$ given by \eqref{eq: diagonal A} where Hypothesis \ref{hyp: 2} holds, satisfy the hard-Lefschetz condition with respect to any symplectic form $\omega$ on $\g$. 
\end{theorem}
\begin{proof}
    Similar to Theorem \ref{thm: main result I}.
\end{proof}

\smallskip

\section{Lattices} \label{sec: lattices} 

In this section we will show the existence of lattices in the Lie groups associated to the Lie algebras we have been working with, for an appropriate choice of the parameters $b_2,\ldots,b_{2n}$. 
The main tool we will use is the following useful result, which provides a criterion to determine when an almost abelian Lie group admits lattices. 

\begin{proposition}\label{prop:latt}\cite{Bo}
Let $G=\R\ltimes_\phi\R^d$ be a unimodular almost abelian Lie group. Then $G$ admits a lattice if and only if there exists  $t_0\neq 0$ such that $\phi(t_0)$ is conjugate to a matrix in $\operatorname{SL}(d,\Z)$.  In this situation, a lattice is given by $\Gamma=t_0 \Z\ltimes P\mathbb Z^{d}$, where $P\in \operatorname{GL}(d,\R)$ satisfies $P^{-1}\phi(t_0)P\in \operatorname{SL}(d,\Z)$. 
\end{proposition}

\begin{remark}
    Note that if $E:=P^{-1}\phi(t_0)P$ then $\Gamma\cong \Z\ltimes_E \Z^d$, where the group multiplication in this last group is given by 
    \[ (m,(p_1,\ldots,p_d))\cdot (n,(q_1,\ldots, q_d))=(m+n,(p_1,\ldots,p_d)+E^m(q_1,\ldots, q_d)). \]
\end{remark}

\medskip

We begin with the Case II, which is easier.

\subsection{Case II} The defining matrix $A\in \operatorname{GL}(2n-1,\R)$ of this case, 
\begin{align*}
    A= \operatorname{diag}(0,\underbrace{1\ldots, 1}_{n-1},\underbrace{-1,\ldots, -1}_{n-1})\in\mathfrak{sl}(2n-1,\R),
\end{align*}
\noindent is conjugate to
\begin{align*}
    A'=\operatorname{diag}(0,1,-1,\ldots, 1,-1)\in\mathfrak{sl}(2n-1,\R),
\end{align*}
\noindent which exponentiates to
\begin{align*}
    \e^{tA'}=\operatorname{diag}(0,\e^t,\e^{-t},\ldots, \e^t,\e^{-t})\in\operatorname{SL}(2n-1,\R).
\end{align*}
\noindent For $m\in \N$ such that $m\geq 3$, the numbers
\begin{align}\label{eq: tm}
    t_m:=\log \frac{m+\sqrt{m^2-4}}{2}, \quad t_m^{-1} := \log \frac{m - \sqrt{m^2-4}}{2}
\end{align}
\noindent exponentiate to the roots of $p_m(x) = x^2 - mx + 1$. Each $2 \times 2$ block $\begin{bmatrix}e^{t_m}&0\\ 0& e^{-t_m} \end{bmatrix}$ in $\e^{t_mA'}$ is then conjugate to the companion matrix $\begin{bmatrix}0&-1\\ 1& m \end{bmatrix}$ of $p_m$, and thus $\e^{t_mA'}$ is conjugate to a matrix in $\operatorname{SL}(2n-1,\Z)$.  The following result is then a consequence of Proposition \ref{prop:latt}, Hattori's theorem (see equation \eqref{deRham}), and Theorem \ref{thm: main result I}. 

\begin{theorem}
Let $G_{\mathrm{II}}$ denote the simply connected almost abelian Lie group whose Lie algebra satisfies Hypothesis \ref{hyp: 2}. For any $m\in \N$ such that $m\geq 3$ there is a lattice $\Gamma_m$ in $G_{ \mathrm{II} }$. The solvmanifolds $\Gamma_m \backslash G_{\mathrm{II}}$, equipped with any invariant symplectic form, satisfy the hard-Lefschetz condition. Moreover, the Betti numbers of $\Gamma_m\backslash G_{\mathrm{II}}$ are independent of $m$, and are given by the values in Corollary \ref{cor: betti numbers of case II}.
\end{theorem}

\smallskip

\subsection{Case I} The idea in this case is to show that, for any $n\geq 2$, we can choose $n-1$ values $m_2,\ldots, m_{n}\in\N$, $m_j\geq 3$, such that the real numbers $t_{m_2},\ldots, t_{m_n}$ given by \eqref{eq: tm} satisfy Hypothesis \ref{hyp: 1}.

Given $m\geq 3$, let us denote $s_m:=\frac{m+\sqrt{m^2-4}}{2}$, so that $s_m=\e^{t_m}$. Moreover, let us write $m^2-4=k_m^2 d_m$, where $k_m\in \N$ and $d_m\in \N$ is square-free. As a consequence we have that $s_m\in \Q(\sqrt{d_m})$.

We will show next that we can choose $m_2,\ldots, m_n$ such that the corresponding square-free numbers $d_{m_2}, \ldots, d_{m_n}$ are pairwise coprime. In fact, we will begin the other way round: let $d_{2}, \ldots, d_{n}$ be pairwise coprime square-free natural numbers, with $d_{j}\geq 2$ for all $j$. Consider now the following generalized Pell's equations:
\begin{equation}\label{eq: Pell}
x^2-d_{j}y^2=4, \quad  j=2,\ldots,n. 
\end{equation}
It follows from \cite[Section 4.4.2]{AA} that the equations \eqref{eq: Pell} have integer solutions for any $d_j$. Let $(m_j,k_j)\in\N \times \N$ be a solution of $\eqref{eq: Pell}$, different from the trivial one $(x,y)=(2,0)$; therefore we have that $m_j\geq 3$ and $k_j\geq 1$. These numbers $m_2,\ldots,m_n$ are the ones that we need.

\begin{lemma} \label{lemma: estos sirven}
    If $m_2, \ldots, m_n$ are as above then the numbers $t_{m_2}, \ldots, t_{m_n}$ satisfy Hypothesis \ref{hyp: 1}. 
\end{lemma} 
\begin{proof}
Let us assume that there exists a linear combination 
\begin{equation}\label{eq: sum}
  \sum_{j=2}^n \varepsilon_jt_{m_j}=0  
\end{equation}
with $\varepsilon_j=\pm 1$ for all $j$; we will arrive at a contradiction. Applying the exponential function to both sides of this equality we get 
\begin{equation}\label{eq: prod}
    \prod_{j=2}^n s_{m_j}^{\varepsilon_j}=1.
\end{equation} 
Note that
\[s_{m_j}^{\varepsilon_j}=\frac12\left(m_j+\varepsilon_j \sqrt{m_j^2-4}\right)=\frac12\left(m_j+\varepsilon_j k_j\sqrt{d_j}\right)\in \Q(\sqrt{d_j}).\]
Let us denote $M:=m_2\cdots m_n$. Expanding the product in \eqref{eq: prod} we obtain
\begin{equation}\label{eq: expand}
     M+\sum_{i=2}^n \frac{M}{m_i}\varepsilon_i k_i\sqrt{d_i} + \sum_{2\leq i<j\leq n} \frac{M}{m_i m_j}\varepsilon_i\varepsilon_j k_i k_j\sqrt{d_i d_j}+\cdots +\varepsilon_2\cdots \varepsilon_n k_2\cdots k_n \sqrt{d_2\cdots d_n}=2^{n-1}.
\end{equation}
Note that the integer number $2^{n-1}-M$ is non-zero, since $M=m_2\cdots m_n\geq 3^{n-1}$. Since the numbers $d_2,\ldots, d_n$ are square-free and pairwise coprime, we have that all the radicals in \eqref{eq: expand} are different and belong to the set $D:=\{ \sqrt{d} \mid d\in \N \text{ is square-free} \}$. However, it is well known\footnote{See for instance the discussion in https://qchu.wordpress.com/2009/07/02/square-roots-have-no-unexpected-linear-relationships/ and the references therein.} that no non-zero integer can be represented as a nontrivial linear combination of radicals in $D$, and this is in contradiction with \eqref{eq: expand}.
Hence, the only possibility for \eqref{eq: sum} to hold is $\varepsilon_j=0$ for all $j$, and this concludes the proof.
\end{proof}

\medskip

Now, the matrix $A:=\operatorname{diag}(0,t_{m_2},\ldots,t_{m_n},-t_{m_2},\ldots,-t_{m_n})$ is certainly conjugate to  \\  $A':=\operatorname{diag}(0,t_{m_2},-t_{m_2},\ldots,t_{m_n},-t_{m_n})$. Each $2 \times 2$ block $\begin{bmatrix}e^{t_{m_j}}&0\\ 0& e^{-t_{m_j}} \end{bmatrix}$ in the exponential $\e^{A'}$ is then conjugate to $\begin{bmatrix}0&-1\\ 1& m_j \end{bmatrix}$, and thus $\e^{A'}$ is conjugate to a matrix in $\operatorname{SL}(2n-1,\Z)$. The following result is then a consequence of Proposition \ref{prop:latt}, Hattori's theorem (see equation \eqref{deRham}), and Theorem \ref{thm: main result II}. 

\begin{theorem}
Let $G_{\mathrm{I}}$ denote the simply connected almost abelian Lie group whose Lie algebra satisfies Hypothesis \ref{hyp: 1}.    For any $n \in \N$ such that $n \geq 2$ and $m_2, \ldots, m_n$ such as in Lemma \ref{lemma: estos sirven}, denoted collectively as $m := (m_2, \ldots, m_n)$, there is a lattice $\Gamma_{n,m}$ in $G_{ \mathrm{I} }$. The solvmanifolds $\Gamma \backslash G_{\mathrm{I}}$, equipped with any invariant symplectic form, satisfy the hard-Lefschetz condition. Moreover, the Betti numbers of $\Gamma_{n,m} \backslash G_{\mathrm{I}}$ are independent of $n$ and $m$, and are given by the values in Corollary \ref{cor: betti numbers of case I}.
\end{theorem}

\begin{example}
Let us set $n=5$ and $d_2=2,\,d_3=3,\,d_4=5,\, d_5=7$. We may choose as solutions of the associated Pell's equations \eqref{eq: Pell} the following pairs:
\[ (m_2,k_2)=(6,4),\, (m_3,k_3)=(4,2),\, (m_4,k_4)=(3,1), \,(m_5,k_5)=(16,6).\]
These values of $m_2,\ldots,m_5$ determine a 10-dimensional almost abelian Lie group $G_{\mathrm{I}}$ satisfying Hypothesis \ref{hyp: 1}, and $G_{\mathrm{I}}$ has a lattice $\Gamma_{\mathrm{I}}$. This lattice is isomorphic to $\Z\ltimes_E \Z^9$, where $E\in \operatorname{SL}(9,\Z)$ is given by
\[ E=[1]\oplus \begin{bmatrix}
    0&-1\\ 1&6
\end{bmatrix}\oplus \begin{bmatrix}
    0&-1\\1&4
\end{bmatrix}\oplus \begin{bmatrix}
    0&-1\\1&3
\end{bmatrix}\oplus \begin{bmatrix}
    0&-1\\1&16
\end{bmatrix}.\]
\end{example}

\

\begin{remark}
There is another way in which we can show, for any $n\geq 2$, the existence of numbers $m_2,\ldots, m_{n}\in\N$ such that the real numbers $t_{m_2},\ldots, t_{m_n}$ given by \eqref{eq: tm} satisfy Hypothesis \ref{hyp: 1}. We describe it as follows. 

Note first that $t_m=\cosh^{-1}(m/2)$. Consider, for $n\geq 2$, an increasing sequence of natural numbers $1\leq k_2<k_3<\cdots <k_n$, and let us choose, for $2\leq j \leq n$, natural numbers $m_j$ satisfying
\[ 2\cosh(2^{k_j-1})\leq m_j < 2\cosh(2^{k_j}).\]
Then the real numbers $t_{m_j}$ satisfy $2^{k_j-1}<t_{m_j}<2^{k_j}$ for $j\geq 2$, and therefore the binary representation of $\floor{t_{m_j}}$ has $k_j$ digits. It follows that $\{t_{m_j}\,:\, j=2,\ldots,n\}$ satisfy Hypothesis \ref{hyp: 1}. The construction of the associated almost abelian Lie group and its lattice is exactly the same as before. 

For instance, if $n=5$ and setting $k_j=j-1$ for $j=2,\ldots, 5$, we may choose 
\[ m_2=4, \quad m_3=8, \quad m_4=55, \quad m_5=2981.  \]
Even though this choice of $m_2, \ldots,m_n$ guarantees the existence of lattices for any $n$, the disadvantage of this method is that the numbers $m_n$ grow extremely fast with $n$, so that it is unfeasible to write them down explicitly.
\end{remark}

\smallskip

\textit{Declaration of competing interest.} 
The authors declare that they have no known competing financial interests or personal relationships that could have appeared to influence the work reported in this paper.

\smallskip

\textit{Acknowledgements.} This work was partially supported by CONICET, SECyT-UNC and ANPCyT (Argentina). The authors are grateful to Facundo Javier Gelatti, Emilio Lauret, Juan Pablo Rossetti and Diego Sulca for useful comments and suggestions.

\

\end{document}